\documentclass[a4paper]{amsart}

\usepackage[german,english]{babel}
\usepackage{graphicx}
\RequirePackage{amsmath}

\RequirePackage{amsmath}
\RequirePackage{bm}
\RequirePackage{amssymb}
\RequirePackage{upref}
\RequirePackage{amsthm}
\RequirePackage{enumerate}
\RequirePackage{pb-diagram}
\RequirePackage{amsfonts}
\RequirePackage[mathscr]{eucal}
\RequirePackage{verbatim}
\RequirePackage{xr}
\RequirePackage{graphicx}
\usepackage{calc}
\usepackage{xspace}

\renewcommand{\qedhere}{\relax}

\newcommand{\cf}{cf.\@\xspace}
\newcommand{\resp}{resp.\@\xspace}

\newcommand{\aev}{a.e.\@\xspace}

\newcommand{\al}{\alpha}
\newcommand{\bet}{\beta}
\newcommand{\ga}{\gamma}
\newcommand{\de}{\delta }
\newcommand{\e}{\epsilon}

\newcommand{\f}{\varphi}
\newcommand{\h}{\eta}

\newcommand{\ka}{\kappa}
\newcommand{\lam}{\lambda}
\newcommand{\m}{\mu}
\newcommand{\n}{\nu}

\newcommand{\s}{\sigma}
\newcommand{\x}{\xi}

\newcommand{\C}{\varGamma}
\newcommand{\D}{\varDelta}
\newcommand{\F}{\varPhi}
\newcommand{\Lam}{\varLambda}
\newcommand{\Om}{\varOmega}

\newcommand{\di}[1]{#1\nobreakdash-\hspace{0pt}dimensional}
\newcommand{\nbdd}{\nobreakdash--}

\newcommand{\fv}[2]{#1\hspace{0pt}_{|_{#2}}}

\newcommand{\so}{{\mc S_0}}

\newcommand{\const}{\tup{const}}
\newcommand{\ndash}{\nobreakdash--}

\newcommand{\msp[1]}[1]{\mspace{#1mu}}

\newcommand{\R}[1][n+1]{{\protect\mathbb R}^{#1}}

\newcommand{\N}{{\protect\mathbb N}}

\newcommand{\eR}{\stackrel{\lower1ex \hbox{\rule{6.5pt}{0.5pt}}}{\msp[3]\R[]}}
\newcommand{\eN}{\stackrel{\lower1ex \hbox{\rule{6.5pt}{0.5pt}}}{\msp[1]\N}}
\newcommand{\eO}{\stackrel{\lower1ex
\hbox{\rule{6pt}{0.5pt}}}{\msc O}}

\DeclareMathOperator{\diam}{diam}

\DeclareMathOperator{\graph}{graph}

\newcommand\im{\implies}
\newcommand\ra{\rightarrow}

\newcommand{\ua}{\uparrow}
\newcommand{\da}{\downarrow}

\newcommand\pa{\partial}
\newcommand\pde[2]{\frac {\partial#1}{\partial#2}}
 
\newcommand\sql[1][u]{\sqrt{1-|D#1|^2}}
\newcommand{\un}{\infty}
\newcommand{\A}{\forall}

\newcommand{\set}[2]{\{\,#1\colon #2\,\}}
\newcommand{\uu}{\cup}
\newcommand{\ii}{\cap}
\newcommand{\uuu}{\bigcup}

\newcommand{\uud}{ \stackrel{\lower 1ex \hbox {.}}{\uu}}
\newcommand{\uuud}[1]{ \stackrel{\lower 1ex \hbox {.}}{\uuu_{#1}}}
\newcommand\su{\subset}

\newcommand\eS{\emptyset}
\newcommand{\sminus}[1][28]{\raise 0.#1ex\hbox{$\scriptstyle\setminus$}}

\newcommand{\wed}{\wedge}

\newcommand{\abs}[1]{\lvert#1\rvert}

\newcommand{\norm}[1]{\lVert#1\rVert}

\newcommand{\nnorm}[1]{| \mspace{-2mu} |\mspace{-2mu}|#1| \mspace{-2mu}
|\mspace{-2mu}|}
\newcommand{\spd}[2]{\protect\langle #1,#2\protect\rangle}

\newcommand\ch[3]{\varGamma_{#1#2}^#3}
\newcommand\cha[3]{{\bar\varGamma}_{#1#2}^#3}

\newcommand{\riem}[4]{R_{#1#2#3#4}}
\newcommand{\riema}[4]{{\bar R}_{#1#2#3#4}}

\newcommand{\tit}{\textit}

\newcommand{\tup}{\textup}

\newcommand{\mc}{\protect\mathcal}
\newcommand{\msc}{\protect\mathscr}

\providecommand{\bysame}{\makebox[3em]{\hrulefill}\thinspace}

\newcommand{\ci}{\cite}

\newcommand{\bt}{\begin{thm}}
\newcommand{\bl}{\begin{lem}}
\newcommand{\bc}{\begin{cor}}
\newcommand{\bd}{\begin{definition}}
\newcommand{\bpp}{\begin{prop}}
\newcommand{\br}{\begin{rem}}
\newcommand{\bn}{\begin{note}}
\newcommand{\be}{\begin{ex}}
\newcommand{\bes}{\begin{exs}}
\newcommand{\bb}{\begin{example}}
\newcommand{\bbs}{\begin{examples}}
\newcommand{\ba}{\begin{axiom}}
\newcommand{\bas}{\begin{assumption}}

\newcommand{\et}{\end{thm}}
\newcommand{\el}{\end{lem}}
\newcommand{\ec}{\end{cor}}
\newcommand{\ed}{\end{definition}}
\newcommand{\epp}{\end{prop}}
\newcommand{\er}{\end{rem}}
\newcommand{\en}{\end{note}}
\newcommand{\ee}{\end{ex}}
\newcommand{\ees}{\end{exs}}
\newcommand{\eb}{\end{example}}
\newcommand{\ebs}{\end{examples}}
\newcommand{\ea}{\end{axiom}}
\newcommand{\eas}{\end{assumption}}

\newcommand{\bp}{\begin{proof}}
\newcommand{\ep}{\end{proof}}
\newcommand{\eps}{\renewcommand{\qed}{}\end{proof}}

\newcommand{\bal}{\begin{align}}

\newcommand{\bi}[1][1.]{\begin{enumerate}[\upshape #1]}
\newcommand{\bia}[1][(1)]{\begin{enumerate}[\upshape #1]}
\newcommand{\bin}[1][1]{\begin{enumerate}[\upshape\bfseries #1]}
\newcommand{\bir}[1][(i)]{\begin{enumerate}[\upshape #1]}
\newcommand{\bic}[1][(i)]{\begin{enumerate}[\upshape\hspace{2\cma}#1]}
\newcommand{\bis}[2][1.]{\begin{enumerate}[\upshape\hspace{#2\parindent}#1]}
\newcommand{\ei}{\end{enumerate}}

\newcommand\ndots{\raise 0.47ex \hbox {,}\hskip0.06em\cdots %
     \raise 0.47ex \hbox {,}\hskip0.06em} 

\newcommand{\q}{\quad}
\newcommand{\qq}{\qquad}

\newcommand{\hp}{\hphantom}

\newcommand\nd{\noindent}

\newskip\Csmallskipamount                                                
\Csmallskipamount=\smallskipamount
\newskip\Cmedskipamount
\Cmedskipamount=\medskipamount
\newskip\Cbigskipamount
\Cbigskipamount=\bigskipamount

\newcommand\cvs{\vspace\Csmallskipamount}   
\newcommand\cvm{\vspace\Cmedskipamount}

\newskip\csa
\csa=\smallskipamount

\newskip\cma
\cma=\medskipamount

\newskip\cba
\cba=\bigskipamount

\newdimen\spt
\newcommand\citem{\cvs\advance\itemno by
1{(\romannumeral\the\itemno})\hskip3pt}
\newcommand{\bitem}{\cvm\nd\advance\itemno by
1{\bf\the\itemno}\hspace{\cma}}

\newcount\itemno
\newcommand{\las}[1]{\label{S:#1}}

\newcommand{\lae}[1]{\label{E:#1}}
\newcommand{\lat}[1]{\label{T:#1}}
\newcommand{\lal}[1]{\label{L:#1}}
\newcommand{\lad}[1]{\label{D:#1}}
\newcommand{\lac}[1]{\label{C:#1}}
\newcommand{\lan}[1]{\label{N:#1}}

\newcommand{\lar}[1]{\label{R:#1}}

\newcommand{\laas}[1]{\label{Ass:#1}}

\newcommand{\rs}[1]{Section~\ref{S:#1}}

\newcommand{\rt}[1]{Theorem~\ref{T:#1}}
\newcommand{\rl}[1]{Lemma~\ref{L:#1}}
\newcommand{\rd}[1]{Definition~\ref{D:#1}}
\newcommand{\ras}[1]{Assumption~\ref{Ass:#1}}
\newcommand{\rc}[1]{Corollary~\ref{C:#1}}

\newcommand{\re}[1]{\eqref{E:#1}}
\newcommand{\frt}[1]{Theorem~\ref{T:#1} on page~\tup{\pageref{T:#1}}}
\newcommand{\frl}[1]{Lemma~\ref{L:#1} on page~\tup{\pageref{L:#1}}}
\newcommand{\frr}[1]{Remark~\ref{R:#1} on page~\tup{\pageref{R:#1}}}
\newcommand{\frd}[1]{Definition~\ref{D:#1} on page~\tup{\pageref{D:#1}}}

\newskip\thmskip
\thmskip=\parindent

\newskip\hsk
\newenvironment{hinw}{\labelsep=0pt\begin{list}{}{\labelsep=0pt\itemindent=0pt\labelwidth=0pt\leftmargin=\parindent\rightmargin=0pt\partopsep=\cba}%
\item\it\nopagebreak\nopagebreak}%
{\end{list}}

\newcommand\bh{\begin{hinw}}
\newcommand{\eh}{\end{hinw}}

\newtheoremstyle{normal}
  {\cba}
  {\cba}
  {}
  {\thmskip}
  {\bfseries}
  {.}
  {\hsk}
  {}

\newtheoremstyle{abschnitt}
  {\cba}
  {\cba}
  {}
  {\thmskip}
  {\bfseries}
  {.}
  {\hsk}
  {}

\newtheoremstyle{italic}
  {\cba}
  {\cba}
  {\itshape}
  {\thmskip}
  {\bfseries}
  {.}
  {\hsk}
  {}

\newtheoremstyle{aufgaben}
  {\cba}
  {\cba}
  {}
  {}
  {\normalsize\bfseries}
  {.}
  {\hsk}
  {}

\newtheoremstyle{break}
  {\cba}
  {\cba}
  {\itshape}
  {}
  {\bfseries}
  {.}
  {\newline}
  {}

\swapnumbers
\theoremstyle{italic}
\newtheorem{thm}[subsection]{Theorem}
\newtheorem{lem}[subsection]{Lemma}
\newtheorem{prop}[subsection]{Proposition}
\newtheorem{cor}[subsection]{Corollary}

\theoremstyle{normal}
\newtheorem{rem}[subsection]{Remark}
\newtheorem{definition}[subsection]{Definition}
\newtheorem{example}[subsection]{Example}
\newtheorem{examples}[subsection]{Examples}
\newtheorem{ex}[subsection]{Exercise}
\newtheorem{note}[subsection]{}
\newtheorem{axiom}[subsection]{Axiom}
\newtheorem{assumption}[subsection]{Assumption}

\theoremstyle{aufgaben}
\newtheorem{exs}[subsection]{Exercises}

\swapnumbers

\numberwithin{equation}{section}
\numberwithin{figure}{section}

\newenvironment{textequation}[1][0.8]
{\begin{equation}
\begin{aligned}
\begin{minipage}{#1\linewidth}}
{\end{minipage}
\end{aligned}
\end{equation}
\ignorespacesafterend}

\newcommand{\btext}{\begin{textequation}}
\newcommand{\etext}{\end{textequation}}
\makeatletter
\def\hinweis{\@startsection{subsection}{2}%
 \z@{0.7\linespacing\@plus 0.5\linespacing}{0.7\linespacing}%
{\normalfont\itshape\indent}}
\makeatother
\RequirePackage{bm}
\RequirePackage{amssymb}
\RequirePackage{upref}
\RequirePackage{amsthm}
\RequirePackage{enumerate}
\RequirePackage{pb-diagram}
\RequirePackage{amsfonts}
\RequirePackage[mathscr]{eucal}
\RequirePackage{verbatim}
\RequirePackage{xr}
\RequirePackage{graphicx}
\usepackage{calc}
\usepackage{xspace}

\newcommand{\ind}[1]{#1}
\newcommand{\indexs}[1]{\relax}
\renewcommand{\index}[1]{\relax}
\newcommand{\inds}[1]{#1}
\newcommand{\fre}[1]{\eqref{E:#1} on page~\tup{\tup{\pageref{E:#1}}}}
\newcommand{\lab}[1]{\label{B:#1}}

\newcommand{\frs}[1]{Section~\ref{S:#1} on page~\tup{\pageref{S:#1}}}

\newcommand{\frn}[1]{Note~\ref{N:#1} on page~\tup{\pageref{N:#1}}}
\newcommand{\rb}[1]{Example~\ref{B:#1}}

\RequirePackage{upref}
\RequirePackage{amsthm}
\RequirePackage{enumerate}
\usepackage[mathscr]{eucal}

\usepackage{xr-hyper}

\usepackage{calc}

\newlength{\oddsidemarginlength}
\newlength{\topmarginlength}

\newcounter{numberoflines}
\newcounter{tempcc}
\oddsidemargin=\oddsidemarginlength
\evensidemargin=\oddsidemargin
\topmargin=\topmarginlength
\begin{document}

\flushbottom


\title{Curvature flows and CMC hypersurfaces}

\author{Claus Gerhardt}
\address{Ruprecht-Karls-Universit\"at, Institut f\"ur Angewandte Mathematik,
Im Neuenheimer Feld 294, 69120 Heidelberg, Germany}
\email{gerhardt@math.uni-heidelberg.de}
\urladdr{http://www.math.uni-heidelberg.de/studinfo/gerhardt/}
\thanks{This research was supported by the Deutsche Forschungsgemeinschaft.}

%
\subjclass[2000]{35J60, 53C21, 53C44, 53C50, 58J05}
\keywords{curvature flows, CMC hypersurfaces, foliations, prescribed curvature, semi-Riemannian manifold, cosmological spacetime, general relativity}
\date{\today}
%


\begin{abstract}
We give an overview of the existence and regularity results for curvature flows and how these flows can be used to solve some problems in geometry and physics.
\end{abstract}

\maketitle

\tableofcontents

\setcounter{section}{0}
\section{Introduction}
We want to give a survey of the existence and regularity results for extrinsic curvature flows in semi-Riemannian manifolds, i.e., Riemannian or Lorentzian ambient spaces. In order to treat both cases simultaneously terminology like spacelike, timelike, etc., that only makes sense in a Lorentzian setting should be ignored in the Riemannian case.

Those, who are already experts in the field, might still be interested in the regularity result in \rt{6.5.5}---especially the time independent $C^{m+2,\al}$-estimates---for converging curvature flows that are graphs,  and in the general curvature estimates for flows in Riemannian manifolds in \rs{h6}. 

\section{Notations and preliminary results}\las 2

The main objective of this section is to state the equations of Gau{\ss}, Codazzi,
and Weingarten for hypersurfaces. In view of the subtle but important
difference  that is to be seen in the \tit{Gau{\ss} equation} depending on the
nature of the ambient space---Riemannian or Lorentzian---, we shall formulate the governing equations of a
hypersurface $M$ in a semi-Riemannian \di{(n+1)} space $N$, which is either
Riemannian or Lorentzian. Geometric quantities in $N$ will be denoted by
$(\bar g_{\alpha\beta}),(\riema \alpha\beta\gamma\delta)$, etc., and those in $M$ by $(g_{ij}), (\riem
ijkl)$, etc. Greek indices range from $0$ to $n$ and Latin from $1$ to $n$; the
summation convention is always used. Generic coordinate systems in $N$ resp.
$M$ will be denoted by $(x^\alpha)$ resp. $(\x^i)$. Covariant differentiation will
simply be indicated by indices, only in case of possible ambiguity they will be
preceded by a semicolon, i.e., for a function $u$ in $N$, $(u_\alpha)$ will be the
gradient and
$(u_{\alpha\beta})$ the Hessian, but e.g., the covariant derivative of the curvature
tensor will be abbreviated by $\riema \alpha\beta\gamma{\delta;\e}$. We also point out that
\begin{equation}
\riema \alpha\beta\gamma{\delta;i}=\riema \alpha\beta\gamma{\delta;\e}x_i^\e
\end{equation}
with obvious generalizations to other quantities.

Let $M$ be a \tit{spacelike} hypersurface, i.e., the induced metric is Riemannian,
with a differentiable normal $\n$. We define the signature of $\n$, $\sigma=\s(\n)$, by
\begin{equation}
\sigma=\bar g_{\alpha\beta}\n^\alpha\n^\beta=\spd \n\n.
\end{equation}
In case $N$ is Lorentzian, $\sigma=-1$, and $\n$ is timelike.

In local coordinates, $(x^\alpha)$ and $(\x^i)$, the geometric quantities of the
spacelike hypersurface $M$ are connected through the following equations
\begin{equation}\lae{2.3}
x_{ij}^\alpha=-\sigma h_{ij}\n^\alpha
\end{equation}
the so-called \tit{Gau{\ss} formula}. Here, and also in the sequel, a covariant
derivative is always a \tit{full} tensor, i.e.,

\begin{equation}
x_{ij}^\alpha=x_{,ij}^\alpha-\ch ijk x_k^\alpha+\cha \beta\gamma\alpha x_i^\beta x_j^\gamma.
\end{equation}
The comma indicates ordinary partial derivatives.

In this implicit definition the \tit{second fundamental form} $(h_{ij})$ is taken
with respect to $-\sigma\n$.

The second equation is the \tit{Weingarten equation}
\begin{equation}
\n_i^\alpha=h_i^k x_k^\alpha,
\end{equation}
where we remember that $\n_i^\alpha$ is a full tensor.

Finally, we have the \tit{Codazzi equation}
\begin{equation}
h_{ij;k}-h_{ik;j}=\riema\alpha\beta\gamma\delta\n^\alpha x_i^\beta x_j^\gamma x_k^\delta
\end{equation}
and the \tit{Gau{\ss} equation}
\begin{equation}
\riem ijkl=\sigma \{h_{ik}h_{jl}-h_{il}h_{jk}\} + \riema \alpha\beta\gamma\delta x_i^\alpha x_j^\beta x_k^\gamma
x_l^\delta.
\end{equation}
Here, the signature of $\n$ comes into play.

\bd\lad{1.6.7}
(i) Let $F\in C^0(\bar\C)\ii C^{2,\al}(\C)$ be a strictly monotone curvature function, where $\C\su\R[n]$ is a convex, open, symmetric cone containing the positive cone, such that
\begin{equation}
\fv F{\pa \C}=0\q\wed\q \fv F\C>0.
\end{equation}

Let $N$ be semi-Riemannian. A spacelike, orientable\footnote{A hypersurface is said to be orientable, if it has a continuous normal field.} hypersurface $M\su N$ is called \tit{admissible}\index{admissible hypersurface}, if its principal curvatures with respect to a chosen normal lie in $\C$. This definition also applies to subsets of $M$.

\cvm
(ii) Let $M$ be an admissible hypersurface and $f$ a function defined in a neighbourhood of $M$. $M$ is said to be an \tit{\ind{upper barrier}} for the pair $(F,f)$, if
\begin{equation}
\fv FM\ge f
\end{equation}

(iii) Similarly, a spacelike, orientable hypersurface $M$ is called a \tit{\ind{lower barrier}} for the pair $(F,f)$, if at the points $\varSigma\su M$, where
$M$ is admissible, there holds
\begin{equation}
\fv F\varSigma\le f.
\end{equation}
$\varSigma$ may be empty.

\cvm
(iv) If we consider the mean curvature function, $F=H$, then we suppose $F$ to be defined in $\R[n]$ and any spacelike, orientable hypersurface is admissible.
\ed

One of the assumptions that are used when proving a priori estimates is that there exists a strictly convex function 
$\chi\in C^2(\bar \Omega)$ in a given domain $\Om$. We shall state sufficient geometric conditions
guaranteeing the existence of such a function. The lemma below will be valid in Lorentzian as well as Riemannian manifolds, but we formulate and prove it only for the Lorentzian case.

\bl\lal{2.8.1}
Let $N$ be globally hyperbolic, $\so$ a Cauchy hypersurface, $(x^\alpha)$ a special
coordinate system associated with $\so$, and $\bar \Omega\su N$ be compact. Then,
there exists a strictly convex function $\chi\in C^2(\bar \Omega)$ provided the level
hypersurfaces $\{x^0=\const\}$ that intersect $\bar \Omega$ are strictly convex.
\el

\bp
For greater clarity set $t=x^0$, i.e., $t$ is a globally defined time function. Let
$x=x(\x)$ be a local representation for $\{t=\const\}$, and $t_i,t_{ij}$ be the
covariant derivatives of $t$ with respect to the induced metric, and
$t_\alpha,t_{\alpha\beta}$ be the covariant derivatives in $N$, then
\begin{equation}
0=t_{ij}=t_{\alpha\beta}x_i^\alpha x_j^\beta+t_\alpha x_{ij}^\alpha,
\end{equation}
and therefore,
\begin{equation}\lae{2.26}
t_{\alpha\beta}x_i^\alpha x_j^\beta=-t_\alpha x_{ij}^\alpha=-\bar h_{ij} t_\alpha \n^\alpha.
\end{equation}
Here, $(\n^\alpha)$ is past directed, i.e., the right-hand side in \re{2.26} is positive
definite in $\bar \Omega$, since $(t_\alpha)$ is also past directed.

Choose $\lambda>0$ and define $\chi=e^{\lambda t}$, so that
\begin{equation}
\chi_{\alpha\beta}=\lambda^2e^{\lambda t} t_\alpha t_\beta+\lambda e^{\lambda t} t_{\alpha\beta}.
\end{equation}

Let $p\in \Omega$ be arbitrary, $\mc S=\{t=t(p)\}$ be the level hypersurface through
$p$, and $(\h^\alpha)\in T_p(N)$. Then we conclude
\begin{equation}
e^{-\lambda t}\chi_{\alpha\beta}\h^\alpha \h^\beta=\lambda^2\abs{\h^0}^2+\lambda t_{ij}\h^i\h^j+2\lambda t_{0j}\h^0\h^i,
\end{equation}
where $t_{ij}$ now represents the left-hand side in  \re{2.26}, and we infer further
\begin{equation}
\begin{aligned}
e^{-\lambda t}\chi_{\alpha\beta}\h^\alpha\h^\beta&\ge \tfrac{1}{2} \lambda^2 {\abs\h^0}^2 +[\lambda
\e-c_\e]\sigma_{ij}\h^i\h^j\\
&\ge \tfrac{\e}{2}\lambda \{-\abs{\h^0}^2+\sigma_{ij}\h^i\h^j\}
\end{aligned}
\end{equation}
for some $\e>0$, and where $\lambda$ is supposed to be large. Therefore, we
have in $\bar \Omega$
\begin{equation}
\chi_{\alpha\beta}\ge c\bar g_{\alpha\beta}\,,\q c>0,
\end{equation}
i.e., $\chi$ is strictly convex.
\ep

\section{Evolution equations for some geometric quantities}\las{1.3}

Curvature flows are used for different purposes, they can be merely vehicles to approximate a stationary solution, in which case the flow is driven not only by a curvature function but also by the corresponding right-hand side, an external force, if you like, or the flow is a pure curvature flow driven only by a curvature function, and it is used to analyze the topology of the initial hypersurface, if the ambient space is Riemannian, or the singularities of the ambient space, in the Lorentzian case.

In this section we are treating very general curvature flows\footnote{We emphasize that we are only considering flows driven by the extrinsic curvature not by the intrinsic curvature.} in a semi-Riemann\-ian manifold $N=N^{n+1}$, though we only have the Riemannian or Lorentzian case in mind, such that the flow can be either a pure curvature flow or may also be driven by an external force. The nature of the ambient space, i.e., the signature of its metric, is expressed by a parameter $\s=\pm 1$, such that $\s=1$ corresponds to the Riemannian and $\s=-1$ the Lorentzian case. The parameter $\s$ can also be viewed as the signature  of the normal of the spacelike hypersurfaces, namely, 
\begin{equation}
\s=\spd\n\nu.
\end{equation}

 Properties like spacelike, achronal, etc., however, 
only make sense, when $N$ is Lorentzian and should be ignored otherwise.
 
 We consider a strictly monotone, symmetric, and concave curvature $F\in C^{4,\al}(\C)$, homogeneous of degree $1$, a function $0<f\in C^{4,\al}(\Om)$,  where $\Om\su N$ is an open set, and a real function $\F\in C^{4,\al}(\R[]_+)$ satisfying
\begin{equation}
\dot\F>0\q \tup{and}\q \ddot\F\le 0.
\end{equation}

For notational reasons, let us abbreviate
\begin{equation}
\tilde f=\F(f).
\end{equation}

Important examples of functions $\F$ are
 \begin{equation}
\F(r)=r,\q
\F(r)=\log r,\q \F(r)=-r^{-1}
\end{equation}
or
\begin{equation}
\F(r)=r^{\frac1k},\q \F(r)=-r^{-\frac1k},\qq k\ge1.
\end{equation}
\br
The latter choices are necessary, if the curvature function $F$ is not homogeneous of degree $1$ but of degree $k$, like the symmetric polynomials $H_k$. In this case we would sometimes like to define $F=H_k$ and not $H_k^{1/k}$, since 
\begin{equation}
F^{ij}=\frac{\pa F}{\pa h_{ij}}
\end{equation}
is then divergence free, if the ambient space is a spaceform,  though on the other hand we need a concave operator for technical reasons, hence we have to take the $k$-th root.
\er
The curvature flow is given by  the evolution problem
\begin{equation}\lae{1.3.5}
\begin{aligned}
\dot x&=-\s(\F-\tilde f)\n,\\
x(0)&=x_0,
\end{aligned}
\end{equation}
where $x_0$ is an embedding of an initial  compact, spacelike
hypersurface $M_0\su \Om$ of class $C^{6,\al}$, $\F=\F(F)$, and $F$ is evaluated at the principal curvatures
of the flow hypersurfaces $M(t)$, or, equivalently, we may assume that $F$
depends on the second fundamental form $(h_{ij})$ and the metric $(g_{ij})$ of
$M(t)$; $x(t)$ is the embedding of $M(t)$ and $\sigma$ the signature of the  normal $\n=\n(t)$, which is identical to the normal used in the \ind{Gaussian formula} \fre{2.3}.

The initial hypersurface should be \tit{admissible}\index{admissible hypersurface}, i.e., its principal curvatures should belong to the convex, symmetric cone $\C\su \R[n]$.

This is a parabolic problem, so short-time existence is guaranteed, \cf \cite[Chapter 2.5]{cg:cp}

There will be a slight ambiguity in the terminology, since we shall call the
evolution parameter \tit{time}, but this lapse shouldn't cause any
misunderstandings, if the ambient space is Lorentzian.

At the moment we consider a sufficiently smooth solution of the initial value problem \re{1.3.5} and want to show how the metric, the second fundamental form, and the
normal vector of the hypersurfaces $M(t)$ evolve. All time derivatives are
\tit{total} derivatives, i.e., covariant derivatives of tensor fields defined over the curve $x(t)$, \cf \cite[Chapter 11.5]{cg:aII}; $t$ is the flow parameter, also referred to as \tit{\ind{time}}, and $(\xi^i)$ are local coordinates of the initial embedding $x_0=x_0(\xi)$ which will also serve as coordinates for the the flow hypersurfaces $M(t)$. The coordinates in $N$ will be labelled $(x^\al)$, $0\le \al\le n$.
\bl[Evolution of the metric]
The metric $g_{ij}$ of $M(t)$ satisfies the evolution equation
\begin{equation}\lae{1.3.6}
\dot g_{ij}=-2\s(\F-\tilde f)h_{ij}.
\end{equation}
\el
\bp
Differentiating
\begin{equation}
g_{ij}=\spd{x_i}{x_j}
\end{equation}
covariantly with respect to $t$ yields
\begin{equation}
\begin{aligned}
\dot g_{ij}&=\spd{\dot x_i}{x_j}+\spd{x_i}{\dot x_j}\\
&=-2 \s (\F-\tilde f)\spd{x_i}{\nu_j}=-2\s (\F-\tilde f)h_{ij},
\end{aligned}
\end{equation}
in view of the Codazzi equations.
\ep

\bl[Evolution of the normal]
The normal vector evolves according to
\begin{equation}\lae{1.3.9}
\dot \n=\nabla_M(\F-\tilde f)=g^{ij}(\F-\tilde f)_i x_j.
\end{equation}
\el

\bp
Since $\nu$ is unit normal vector we have $\dot\nu\in T(M)$. Furthermore, differentiating 
\begin{equation}
0=\spd\nu{x_i}
\end{equation}
with respect to $t$, we deduce
\begin{equation}
\spd{\dot\nu}{x_i}=-\spd{\nu}{\dot x_i}=(\F-\tilde f)_i.\qedhere
\end{equation}
\ep

\bl[Evolution of the second fundamental form]\lal{1.3.3}
The second fundamental form evolves according to
\begin{equation}\lae{1.3.12}
\dot h_i^j=(\F-\tilde f)_i^j+\sigma (\F-\tilde f) h_i^k h_k^j+\sigma (\F-\tilde f) \riema
\alpha\beta\gamma\delta\n^\alpha x_i^\beta \n^\gamma x_k^\delta g^{kj}
\end{equation}
and
\begin{equation}\lae{1.3.13}
\dot h_{ij}=(\F-\tilde f)_{ij}-\sigma (\F-\tilde f) h_i^k h_{kj}+\sigma (\F-\tilde f) \riema
\alpha\beta\gamma\delta\n^\alpha x_i^\beta \n^\gamma x_j^\delta.
\end{equation}
\el

\bp
We use the \tit{\ind{Ricci identities}} to interchange the covariant derivatives of $\nu$ with respect to $t$ and $\xi^i$
\begin{equation}
\begin{aligned}
\tfrac D{dt}(\n^\al_i)&=(\dot\nu^\al)_i-\bar R^\al_{\hp{\al}\bet\ga\de}\nu^\bet x^\ga_i \dot x^\de\\
&= g^{kl}(\F-\tilde f)_{ki}x^\al_l+g^{kl}(\F-\tilde f)_kx^\al_{li}-\bar R^\al_{\hp{\al}\bet\ga\de} \nu^\bet x^\ga_i \dot x^\de.
\end{aligned}
\end{equation}
For the second equality we used \re{1.3.9}.
On the other hand, in view of the Weingarten equation we obtain
\begin{equation}
\tfrac D{dt}(\nu^\al_i)=\tfrac D{dt}(h^k_i x^\al_k)=\dot h^k_ix^\al_k+h^k_i\dot x^\al_k.
\end{equation}

Multiplying the resulting equation with $\bar g_{\al\bet}x^\bet_j$ we conclude
\begin{equation}
\begin{aligned}
\dot h^k_ig_{kj}-\s(\F-\tilde f)h^k_ih_{kj}&=(\F-\tilde f)_{ij}+\s(\F-\tilde f) \riema \al\bet\ga\de \nu^\al x^\bet_i\nu^\ga x^\de_j
\end{aligned}
\end{equation}
or equivalently \re{1.3.12}.

To derive \re{1.3.13}, we differentiate
\begin{equation}
h_{ij}=h^k_ig_{kj}
\end{equation}
with respect to $t$ and use \re{1.3.6}.  
\ep

We emphasize that equation \re{1.3.12} describes the evolution of the second fundamental form more meaningfully than \re{1.3.13}, since the mixed tensor is independent of the metric.

\bl[Evolution of $(\F-\tilde f)$]\lal{1.3.4}
The term $(\F-\tilde f)$ evolves according to the equation
\begin{equation}\lae{1.3.18}
\begin{aligned}
{(\F-\tilde f)}^\prime-\dot\F F^{ij}(\F-\tilde f)_{ij}=&\msp[3]\sigma \dot \F
F^{ij}h_{ik}h_j^k (\F-\tilde f)
 +\sigma\tilde f_\alpha\n^\alpha (\F-\tilde f)\\
&+\sigma\dot\F F^{ij}\riema \alpha\beta\gamma\delta\n^\alpha x_i^\beta \n^\gamma x_j^\delta (\F-\tilde f),
\end{aligned}
\end{equation}
where
\begin{equation}
(\F-\tilde f)^{\prime}=\frac{d}{dt}(\F-\tilde f)
\end{equation}
and
\begin{equation}
\dot\F=\frac{d}{dr}\F(r).
\end{equation}
\el

\bp
When we differentiate $F$ with respect to $t$ we consider $F$ to depend on the mixed tensor $h^j_i$ and conclude
\begin{equation}
(\F-\tilde f)'=\dot\F F^i_j\dot h^j_i-\tilde f_\al \dot x^\al;
\end{equation}
The equation \re{1.3.18} then follows in view of \re{1.3.5} and \re{1.3.12}.
\ep

\br\lar{1.3.5}
The preceding conclusions, except \rl{1.3.4}, remain valid for flows which do not depend on the curvature, i.e., for flows
\begin{equation}
\begin{aligned}
\dot x&=-\s (-f)\nu=\s f\nu,\\
x(0&)=x_0,
\end{aligned}
\end{equation}
where $f=f(x)$ is defined in an open set $\Om$ containing the initial spacelike hypersurface $M_0$. In the preceding equations we only have to set $\F=0$ and $\tilde f=f$.

The evolution equation for the mean curvature then looks like
\begin{equation}
\dot H=-\D f-\s\{\abs A^2+\bar R_{\al\bet}\nu^\al\nu^\bet\}f,
\end{equation}
where the Laplacian is the Laplace operator on the hypersurface $M(t)$. This is exactly the derivative of the mean curvature operator with respect to normal variations as we shall see in a moment.
\er

But first let us consider the following example.
\bb\lab{1.3.6}
Let $(x^\al)$ be a future directed Gaussian coordinate system in $N$, such that  the metric can be expressed in the form
\begin{equation}
d\bar s^2=e^{2\psi}\{\s(dx^0)^2+\s_{ij}dx^idx^j\}.
\end{equation}
Denote by $M(t)$ the coordinate slices $\{x^0=t\}$, then $M(t)$ can be looked at as the flow hypersurfaces of the flow
\begin{equation}\lae{3.27c}
\dot x=-\s(-e^\psi)\bar\nu,
\end{equation}
where we denote the geometric quantities of the slices by $\bar g_{ij}$, $\bar \nu$, $\bar h_{ij}$, etc. 

Here $x$ is the embedding
\begin{equation}
x=x(t,\xi^i)=(t,x^i).
\end{equation}
Notice that, if $N$ is Riemannian, the coordinate system and the normal are always chosen such that $\nu^0>0$, while, if $N$ is Lorentzian, we always pick the past directed normal.

 Hence the mean curvature of the slices evolves according to
 \begin{equation}\lae{1.3.27}
\dot {\bar H}=-\D e^\psi-\s \{\abs{\bar A}^2+\bar R_{\al\bet}\bar\nu^\al\bar\nu^\bet\}e^\psi.
\end{equation}
\eb

We can now derive the linearization of the mean curvature operator of a spacelike hypersurface, compact or non-compact.

\bn\lan{1.3.7}
Let $M_0\su N$ be a spacelike hypersurface of class $C^4$. We first assume that $M_0$ is compact; then there exists a tubular neighbourhood $\mc U$ and a corresponding \tit{normal} Gaussian coordinate system $(x^\al)$ of class $C^3$ such that $\frac {\pa{}}{\pa x^0}$ is normal to $M_0$. 

Let us consider in 
$\mc U$ of $M_0$ spacelike hypersurfaces $M$ that can be written as graphs
over $M_0$, $M=\graph u$, in the corresponding normal Gaussian coordinate system.
Then the mean curvature of $M$ can be expressed as
\begin{equation}\lae{2.4}
H=\{-\D u+\bar H-\s v^{-2}u^iu^j\bar h_{ij}\}v,
\end{equation}
where $\s=\spd\nu\nu$,  and hence, choosing $u=\e\f$, $\f\in C^2(M_0)$, we deduce
\begin{equation}\lae{2.5}
\begin{aligned}
\frac d{d\e}\fv H {\e=0}&=-\D \f +\dot{\bar H}\f\\[\cma]
&= -\D\f -\s(\abs{\bar A}^2+\bar R_{\al\bet}\nu^\al\nu^\bet)\f,
\end{aligned}
\end{equation}
in view of \re{1.3.27}.

The right-hand side is the \ind{derivative of the mean curvature operator} applied to $\f$. 

If $M_0$ is non-compact,  tubular neighbourhoods exist locally and the relation \re{2.5} will be valid for any $\f\in C^2_c(M_0)$ by using a partition of unity.
\en

The preceding linearization can be immediately generalized to a hypersurface $M_0$ solving the equation
\begin{equation}\lae{3.31c}
\fv F{M_0}=f,
\end{equation}
where $f=f(x)$ is defined in a neighbourhood of $M_0$ and $F=F(h_{ij})$ is curvature operator.

\bl\lal{3.9c}
Let $M_0$ be of class $C^{m,\al}$, $m\ge 2$, $0\le\al\le 1$, satisfy \re{3.31c}. Let $\mc U$ be a (local) tubular neighbourhood of $M_0$, then the linearization of the operator $F-f$ expressed in the normal Gaussian coordinate system $(x^\al)$ corresponding to $\mc U$ and evaluated at $M_0$ has the form
\begin{equation}\lae{3.33c}
-F^{ij}u_{ij}-\s\{F^{ij}h_i^kh_{kj}+F^{ij}\riema \al\bet\ga\de \nu^\al x^\bet_i\nu^\ga x^\de_j+f_\al\nu^\al\}u,
\end{equation}
where $u$ is a function defined in $M_0$, and all geometric quantities are those of $M_0$; the derivatives are covariant derivatives with respect to the induced metric of $M_0$. The operator will be self-adjoint, if $F^{ij}$ is divergence free.
\el

\bp
For simplicity assume that $M_0$ is compact, and let $u\in C^2(M_0)$ be fixed. Then  the hypersurfaces
\begin{equation}
M_\e=\graph(\e u)
\end{equation}
stay in the tubular neighbourhood $\mc U$ for small $\e$, $\abs\e<\e_0$, and their second fundamental forms $(h_{ij})$ can be expressed as
\begin{equation}
v^{-1}h_{ij}=-(\e u)_{ij}+\bar h_{ij},
\end{equation}
where $\bar h_{ij}$ is the second fundamental form of the coordinate slices $\{x^0=\const\}$. 

We are interested in
\begin{equation}
\frac d{d\e}\fv {(F-f)} {\e=0}.
\end{equation}

To differentiate $F$ with respect to $\e$ it is best to consider the mixed form $(h^j_i)$ of the second fundamental form to derive
\begin{equation}
\begin{aligned}
\frac d{d\e}(F-f)=F^i_j\dot h^j_i-\pde f{x^0}u=-F^{ij}u_{ij}+F^i_j\dot {\bar h}^j_i u-\pde f{x^0}u,
\end{aligned}
\end{equation}
where the equation is evaluated at $\e=0$ and $\dot {\bar h}^j_i $ is the derivative of $\bar h^j_i$ with respect to $x^0$.

The result then follows from the evolution equation \re{1.3.12} for the flow \re{3.27c}, i.e., we have to replace $(\F-\tilde f)$ in \re{1.3.12} by $-1$.
\ep

\section{Essential parabolic flow equations}\las{1.4}

From \fre{1.3.12} we deduce with the help of the Ricci identities a parabolic equation
for the second fundamental form
\bl\lal{1.4.1}
The mixed tensor $h_i^j$ satisfies the parabolic equation
\begin{equation}\lae{1.4.1}
\begin{aligned}
&\qq\qq\dot h_i^j-\dot\F F^{kl}h_{i;kl}^j=\\[\cma]
&\hp{=}\;\sigma \dot\F F^{kl}h_{rk}h_l^rh_i^j-\sigma\dot\F F
h_{ri}h^{rj}+\sigma (\F-\tilde f) h_i^kh_k^j\\
&\hp{+}-\tilde f_{\alpha\beta} x_i^\alpha x_k^\beta g^{kj}+\sigma \tilde f_\alpha\n^\alpha h_i^j+\dot\F
F^{kl,rs}h_{kl;i}h_{rs;}^{\hphantom{rs;}j}\\
&\hp{=}+\ddot \F F_i F^j+2\dot \F F^{kl}\riema \alpha\beta\gamma\delta x_m^\alpha x_i ^\beta x_k^\gamma
x_r^\delta h_l^m g^{rj}\\
&\hp{=}-\dot\F F^{kl}\riema \alpha\beta\gamma\delta x_m^\alpha x_k ^\beta x_r^\gamma x_l^\delta
h_i^m g^{rj}-\dot\F F^{kl}\riema \alpha\beta\gamma\delta x_m^\alpha x_k ^\beta x_i^\gamma x_l^\delta h^{mj} \\
&\hp{=}+\sigma\dot\F F^{kl}\riema \alpha\beta\gamma\delta\n^\alpha x_k^\beta\n^\gamma x_l^\delta h_i^j-\sigma\dot\F F
\riema \alpha\beta\gamma\delta\n^\alpha x_i^\beta\n^\gamma x_m^\delta g^{mj}\\
&\hp{=}+\sigma (\F-\tilde f)\riema \alpha\beta\gamma\delta\n^\alpha x_i^\beta\n^\gamma x_m^\delta g^{mj}\\
&\hp{=}+\dot\F F^{kl}\bar R_{\alpha\beta\gamma\delta;\e}\{\n^\alpha x_k^\beta x_l^\gamma x_i^\delta
x_m^\e g^{mj}+\n^\alpha x_i^\beta x_k^\gamma x_m^\delta x_l^\e g^{mj}\}.
\end{aligned}
\end{equation}
\el
\bp
We start with equation \fre{1.3.12} and shall evaluate the term 
\begin{equation}
(\F-\tilde f)^j_i;
\end{equation}
since we are only working with covariant spatial derivatives in the subsequent proof, we may---and shall---consider the covariant form of the tensor
\begin{equation}
(\F-\tilde f)_{ij}.
\end{equation}
First we have
\begin{equation}
\F_i=\dot\F F_i=\dot\F F^{kl}h_{kl;i}
\end{equation}
and
\begin{equation}
\F_{ij}=\dot\F F^{kl}h_{kl;ij}+\Ddot\F F^{kl}h_{kl;i}F^{rs}h_{rs;j}+\dot\F F^{kl,rs}h_{kl,;i}h_{rs;j}.
\end{equation}

Next, we want to replace  $h_{kl;ij}$ by $h_{ij;kl}$. Differentiating the Codazzi equation
\begin{equation}
h_{kl;i}=h_{ik;l}+\bar R_{\al\bet\ga\de}\nu^\al x^\bet_kx^\ga_lx^\de_i,
\end{equation}
where we also used the symmetry of $h_{ik}$, yields
\begin{equation}
\begin{aligned}
&h_{kl;ij}=h_{ik;lj}+\bar R_{\al\bet\ga\de ;\e}\nu^\al x^\bet_kx^\ga_l x^\de_ix^\e_j\\
&+\bar R_{\al\bet\ga\de} \{\nu^\al_j x^\bet_kx^\ga_lx^\de_i+\nu^\al x^\bet_{kj}x^\ga_lx^\de_i+ \nu^\al x^\bet_kx^\ga_{lj}x^\de_i+\nu^\al x^\bet_kx^\ga_lx^\de_{ij}\}.
\end{aligned}
\end{equation}

To replace $h_{kl;ij}$ by $h_{ij;kl}$ we use the Ricci identities
\begin{equation}
h_{ik;lj}=h_{ik;jl}+h_{ak}R^a_{\hp{a}ilj}+h_{ai}R^a_{\hp{a}klj}
\end{equation}
and differentiate once again the Codazzi equation
\begin{equation}
h_{ik;j}=h_{ij;k}+\bar R_{\al\bet\ga\de}\nu^\al x^\bet_ix^\ga_kx^\de_j.
\end{equation}

To replace $\tilde f_{ij}$ we use the chain rule
\begin{equation}
\begin{aligned}
\tilde f_i&=\tilde f_\al x^\al_i,\\
\tilde f_{ij}&=\tilde f_{\al\bet}x^\al_ix^\bet_j+\tilde f_\al x^\al_{ij}.
\end{aligned}
\end{equation}

Then, because of the Gau{\ss} equation, Gaussian formula,  and Weingarten equation,  the symmetry properties of the Riemann curvature tensor and the assumed homogeneity of $F$, i.e.,
\begin{equation}
F=F^{kl}h_{kl},
\end{equation}
we deduce \re{1.4.1} from \fre{1.3.12} after reverting to the mixed representation.
\ep

\br
If we had assumed $F$ to be homogeneous of degree $d_0$ instead of 1, then we
would have to replace the explicit term $F$---occurring twice in the preceding
lemma---by $d_0F$.
\er

If the ambient semi-Riemannian manifold is a space of constant curvature, then the evolution equation of the second fundamental form simplifies considerably, as can be easily verified.

\bl\lal{1.4.3}
Let $N$ be a space of constant curvature $K_N$, then the second fundamental form of the curvature flow \fre{1.3.5} satisfies the  parabolic equation 
\begin{equation}\lae{1.4.12}
\begin{aligned}
\dot h_i^j-\dot\F F^{kl}h_{i;kl}^j&=\s\dot\F F^{kl}h_{rk}h_l^rh_i^j-\s\dot\F F
h_{ri}h^{rj}+\s (\F-\tilde f) h_i^kh_k^j\\
&\hp{=}\;-\tilde f_{\alpha\beta} x_i^\alpha x_k^\beta g^{kj}+\s\tilde f_\alpha\n^\alpha h_i^j+\dot\F
F^{kl,rs}h_{kl;i}h_{rs;}^{\hphantom{rs;}j}\\
&\hp{=}\;+\Ddot\F F_iF^j\\
&\hp{=}\;+K_N\{(\F-\tilde f)\de^j_i+\dot\F F\de^j_i-\dot\F F^{kl}g_{kl}h^j_i\}.
\end{aligned}
\end{equation}
\el

Let us now assume that the open set $\Om\su N$ containing the flow hypersurfaces can be covered by a Gaussian coordinate system $(x^\al)$, i.e., $\Om$ can be topologically viewed as a subset of $I\times \so$, where $\so$ is a compact Riemannian manifold and $I$ an interval. We assume furthermore, that the flow hypersurfaces can be written as graphs over $\so$
\begin{equation}
M(t)=\set{x^0=u(x^i)}{x=(x^i)\in \so};
\end{equation}
we use the symbol $x$ ambiguously by denoting points $p=(x^\al)\in N$ as well as points $p=(x^i)\in \so$ simply by $x$, however, we are careful to avoid confusions.

Suppose that the flow hypersurfaces are given by an embedding $x=x(t,\x)$, where $\x=(\xi^i)$ are local coordinates of a compact manifold $M_0$, which then has to be  homeomorphic to $\so$, then \begin{equation}
\begin{aligned}
x^0&=u(t,\xi)=u(t,x(t,\xi)),\\
x^i&=x^i(t,\xi).
\end{aligned}
\end{equation}

The induced metric can be expressed as
\begin{equation}
\begin{aligned}
g_{ij}=\spd {x_i}{x_j}=\s u_iu_j+\s_{kl}x^k_ix^l_j,
\end{aligned}
\end{equation}
where
\begin{equation}
u_i=u_kx^k_i,
\end{equation}
i.e.,
\begin{equation}
g_{ij}=\{\s u_ku_l+\s_{kl}\}x^k_ix^l_j,
\end{equation}
hence the (time dependent)  \tit{\ind{Jacobian}} $(x^k_i)$ is invertible, and the $(\xi^i)$ can also be viewed as coordinates for $\so$.

Looking at the component $\al=0$ of the flow equation \fre{1.3.5} we obtain a scalar flow equation
\begin{equation}\lae{1.4.18}
\dot u=-e^{-\psi}v^{-1}(\F-\tilde f),
\end{equation}
which is the same in the Lorentzian as well as in the Riemannian case, where
\begin{equation}
v^2=1-\s \s^{ij}u_iu_j,
\end{equation}
and where 
\begin{equation}
\abs{Du}^2=\s^{ij}u_iu_j
\end{equation}
is of course a scalar, i.e., we obtain the same expression regardless, if we use the coordinates $x^i$ or $\xi^i$.

The time derivative in \re{1.4.18} is a total time derivative, if we consider $u$ to depend on $u=u(t,x(t,\xi))$. For the partial time derivative we obtain
\begin{equation}\lae{1.4.21}
\begin{aligned}
\frac {\pa u}{\pa t}&=\dot u-u_k\dot x^k_i\\
&=-e^{-\psi}v(\F-\tilde f),
\end{aligned}
\end{equation}
in view of \fre{1.3.5} and our choice of normal $\n=(\nu^\al)$
\begin{equation}
(\nu^\al)=\s e^{-\psi}v^{-1}(1,-\s u^i),
\end{equation}
where $u^i=\s^{ij}u_j$.

Controlling the $C^1$-norm of the graphs $M(t)$ is tantamount to controlling $v$, if $N$ is Riemannian, and $\tilde v=v^{-1}$, if $N$ is Lorentzian. The evolution equations satisfied by these quantities are also very important, since they are used for the a priori estimates of the second fundamental form.

Let us start with the Lorentzian case.

\bl[Evolution of $\tilde v$]\lal{1.4.4}
Consider the flow \re{1.3.5} in a Lorentzian space $N$ such that the spacelike flow hypersurfaces can be written  as graphs over $\so$. Then, $\tilde v$ satisfies the evolution equation
\begin{equation}\lae{1.4.23}
\begin{aligned}
\dot{\tilde v}-\dot\F F^{ij}\tilde v_{ij}=&-\dot\F F^{ij}h_{ik}h_j^k\tilde v
+[(\F-\tilde f)-\dot\F F]\h_{\alpha\beta}\n^\alpha\n^\beta\\
&-2\dot\F F^{ij}h_j^k x_i^\alpha x_k^\beta \h_{\alpha\beta}-\dot\F F^{ij}\h_{\alpha\beta\gamma}x_i^\beta
x_j^\gamma\n^\alpha\\
&-\dot\F F^{ij}\riema \alpha\beta\gamma\delta\n^\alpha x_i^\beta x_k^\gamma x_j^\delta\h_\e x_l^\e g^{kl}\\
&-\tilde f_\beta x_i^\beta x_k^\alpha \h_\alpha g^{ik},
\end{aligned}
\end{equation}
where $\h$ is the covariant vector field $(\h_\alpha)=e^{\psi}(-1,0,\dotsc,0)$.
\el

\bp
We have $\tilde v=\spd \h\n$. Let $(\x^i)$ be local coordinates for $M(t)$.
Differentiating $\tilde v$ covariantly we deduce
\begin{equation}\lae{1.4.24}
\tilde v_i=\h_{\alpha\beta}x_i^\beta\n^\alpha+\h_\alpha\n_i^\alpha,
\end{equation}
\begin{equation}\lae{1.4.25}
\begin{aligned}
\tilde v_{ij}= &\msp[5]\h_{\alpha\beta\gamma}x_i^\beta x_j^\gamma\n^\alpha+\h_{\alpha\beta}x_{ij}^\beta\n^\alpha\\
&+\h_{\alpha\beta}x_i^\beta\n_j^\alpha+\h_{\alpha\beta}x_j^\beta\n_i^\alpha+\h_\alpha\n_{ij}^\alpha
\end{aligned}
\end{equation}

The time derivative of $\tilde v$ can be expressed as
\begin{equation}\lae{1.4.26}
\begin{aligned}
\dot{\tilde v}&=\h_{\alpha\beta}\msp\dot x^\beta\n^\alpha+\h_\alpha\dot\n^\alpha\\
&=\h_{\alpha\beta}\n^\alpha\n^\beta(\F-\tilde f)+(\F-\tilde f)^k x_k^\alpha\h_\alpha\\
&=\h_{\alpha\beta}\n^\alpha\n^\beta(\F-\tilde f)+\dot\F F^k x_k^\alpha\h_\alpha-{\tilde f}_\beta x_i^\beta
x_k^\alpha g^{ik}\h_\alpha,
\end{aligned}
\end{equation}
where we have used \fre{1.3.9}.

Substituting \re{1.4.25} and \re{1.4.26} in \re{1.4.23}, and simplifying the resulting
equation with the help of the Weingarten and Codazzi equations, we arrive at the
desired conclusion.
\ep

In the Riemannian case we consider a normal Gaussian coordinate system $(x^\al)$, for otherwise we won't obtain a priori estimates for $v$, at least not without additional strong assumptions. We also refer to $x^0=r$ as the \tit{\ind{radial distance function}}.

\bl[Evolution of $v$]\lal{1.4.5}
Consider the flow \re{1.3.5} in a normal Gaussian coordinate system where the $M(t)$ can be written as graphs of a function $u(t)$ over some compact Riemannian manifold $\so$. Then the quantity
\begin{equation}\lae{1.4.28}
v=\sqrt{1+\abs {Du}^2}=(r_\al\nu^\al)^{-1}
\end{equation}
satisfies the evolution equation
\begin{equation}
\begin{aligned}
&\dot v -\dot\F F^{ij}v_{ij}= -\dot\F F^{ij}h_{ik}h^k_jv-2v^{-1}\dot\F F^{ij}v_iv_j\\
&\hp{=}\q +r_{\al\bet}\nu^\al\nu^\bet[(\F-\tilde f)-\dot\F F]v^2+2\dot\F F^{ij}h^k_ir_{\al\bet}x^\al_k x^\bet_j v^2\\
&\hp{=}\q +\dot\F F^{ij}\bar R_{\al\bet\ga\de}\nu^\al x^\bet_i x^\ga_jx^\de_k r_\e x^\e_m g^{mk}v^2\\
&\hp{=}\q +\dot \F F^{ij}r_{\al\bet\ga}\nu^\al x^\bet_i x^\ga _j v^2 +\tilde f_\al x^\al_mg^{mk}r_\bet x^\bet_k v^2.
\end{aligned}
\end{equation}
\el

\bp
Similar to the proof of the previous lemma. 
\ep


The previous problems can be generalized to the case when the right-hand side $f$ is not only defined in $N$ or in $\bar\Om$ but in the \tit{\ind{tangent bundle}} $T(N)$ \resp $T(\bar\Om)$. Notice that the tangent bundle is a manifold of dimension $2(n+1)$, i.e., in a \tit{\ind{local trivialization}} of $T(N)$ $f$ can be expressed in the form
\begin{equation}
f=f(x,\nu)
\end{equation}
with $x\in N$ and $\nu\in T_x(N)$, \cf \cite[Note 12.2.14]{cg:aII}. Thus, the case $f=f(x)$ is included in this general set up. The symbol $\nu$ indicates that in an equation
\begin{equation}\lae{1.4.30}
\fv FM=f(x,\nu)
\end{equation}
we want $f$ to be evaluated at $(x,\nu)$, where $x\in M$ and $\nu$ is the normal of $M$ in $x$.

The \tit{\ind{Minkowski problem}} or \tit{\ind{Minkowski type problems}} are also covered by the present setting, though the Minkowski problem has the additional property that the problem is transformed via the \tit{\ind{Gau{\ss} map}} to a different semi-Riemannian manifold as a \tit{\ind{dual problem}} and solved there. Minkowski type problems have been treated in \cite{cheng-yau:minkowski}, \cite{guan:annals}, \cite{cg:minkowski} and \cite{cg:minkowski2}.
\br
 The equation \re{1.4.30} will be solved by the same methods as in the special case when $f=f(x)$, i.e., we consider the same curvature flow, the evolution equation \fre{1.3.5}, as before.
 
 The resulting evolution equations are identical with the natural exception, that, when $f$ or $\tilde f$ has to be differentiated, the additional argument has to be considered, e.g.,
 \begin{equation}\lae{1.4.31}
\tilde f_i=\tilde f_\al x^\al_i+\tilde f_{\nu^\bet}\nu^\bet_i=\tilde f_\al x^\al_i+\tilde f_{\nu^\bet}x^\bet_k h^k_i
\end{equation}
and
\begin{equation}\lae{1.4.32}
\dot {\tilde f}=\tilde f_\al \dot x^\al+\tilde f_{\nu^\bet}\dot\nu^\bet=-\s(\F-\tilde f)\tilde f_\al\nu^\al+\tilde f_{\nu^\bet}g^{ij}(\F-\tilde f)_i x^\bet_j.
\end{equation}

The most important evolution equations are explicitly stated below. 
\er

Let us first state the evolution equation for $(\F-\tilde f)$.

\bl[Evolution of $(\F-\tilde f)$]\lal{1.4.7}
The term $(\F-\tilde f)$ evolves according to the equation
\begin{equation}\lae{1.4.3.10}
\begin{aligned}
{(\F-\tilde f)}^\prime-\dot\F F^{ij}(\F-\tilde f)_{ij}&=
\s \dot \F
F^{ij}h_{ik}h_j^k (\F-\tilde f)\\
&\hp{=}\;+\s\tilde f_\al\n^\al (\F-\tilde f)
- \tilde f_{\n^\al}x^\al_i(\F- \tilde
f)_jg^{ij}\\
&\hp{=}\;+\s\dot\F F^{ij}\riema \al\bet\ga\de\n^\al x_i^\bet \n^\ga x_j^\de
(\F-\tilde f),
\end{aligned}
\end{equation}where
\begin{equation}
(\F-\tilde f)^{\prime}=\frac{d}{dt}(\F-\tilde f)
\end{equation}
and
\begin{equation}
\dot\F=\frac{d}{dr}\F(r).
\end{equation}
\el

Here is the evolution equation for the second fundamental form.
\bl\lal{1.4.8}
The mixed tensor $h_i^j$ satisfies the parabolic equation
\begin{equation}\lae{1.4.3.13}
\begin{aligned}
\dot h_i^j&-\dot\F F^{kl}h_{i;kl}^j\\[\cma]&=\s \dot\F
F^{kl}h_{rk}h_l^rh_i^j-\s\dot\F F h_{ri}h^{rj}+\s (\F-\tilde f)
h_i^kh_k^j\\[\cma] 
&\hp{+}-\tilde f_{\al\bet} x_i^\al x_k^\bet g^{kj}+\s \tilde f_\al\n^\al
h_i^j-\tilde f_{\al\n^\bet}(x^\al_i x^\bet_kh^{kj}+x^\al_l x^\bet_k h^{k}_i\, g^{lj})\\
&\hp{=}
-\tilde f_{\n^\al\n^\bet}x^\al_lx^\bet_kh^k_ih^{lj}-\tilde f_{\n^\bet} x^\bet_k
h^k_{i;l}\,g^{lj}  +\s\tilde f_{\n^\al}\n^\al h^k_i h^j_k\\
&\hp{=}+\dot\F
F^{kl,rs}h_{kl;i}h_{rs;}^{\hphantom{rs;}j}+2\dot \F F^{kl}\riema \al\bet\ga\de x_m^\al x_i ^\bet x_k^\ga
x_r^\de h_l^m g^{rj}\\
&\hp{=}-\dot\F F^{kl}\riema \al\bet\ga\de x_m^\al x_k ^\bet x_r^\ga x_l^\de
h_i^m g^{rj}-\dot\F F^{kl}\riema \al\bet\ga\de x_m^\al x_k ^\bet x_i^\ga x_l^\de h^{mj} \\
&\hp{=}+\s\dot\F F^{kl}\riema \al\bet\ga\de\n^\al x_k^\bet\n^\ga x_l^\de h_i^j-\s\dot\F F
\riema \al\bet\ga\de\n^\al x_i^\bet\n^\ga x_m^\de g^{mj}\\
&\hp{=}+\s (\F-\tilde f)\riema \al\bet\ga\de\n^\al x_i^\bet\n^\ga x_m^\de g^{mj}+\ddot \F
F_i F^j\\ &\hp{=}+\dot\F F^{kl}\bar R_{\al\bet\ga\de;\e}\{\n^\al x_k^\bet x_l^\ga x_i^\de
x_m^\e g^{mj}+\n^\al x_i^\bet x_k^\ga x_m^\de x_l^\e g^{mj}\}.
\end{aligned}
\end{equation}
\el

The proof is identical to that of  \rl{1.4.1};
we only have to keep in mind that $f$ now also depends on the normal.

If we had assumed $F$ to be homogeneous of degree $d_0$ instead of $1$, then, we
would have to replace the explicit term $F$---occurring twice in the preceding
lemma---by $d_0F$.

\bl[Evolution of $\tilde v$]\lal{1.4.9}
Consider the flow \re{1.3.5} in a Lorentzian space $N$ such that the spacelike flow hypersurfaces can be written  as graphs over $\so$. Then, $\tilde v$ satisfies the evolution equation
\begin{equation}\lae{1.4.4.3}
\begin{aligned}
\dot{\tilde v}-\dot\F F^{ij}\tilde v_{ij}=&-\dot\F F^{ij}h_{ik}h_j^k\tilde v
+[(\F-\tilde f)-\dot\F F]\h_{\al\bet}\n^\al\n^\bet\\
&-2\dot\F F^{ij}h_j^k x_i^\al x_k^\bet \h_{\al\bet}-\dot\F F^{ij}\h_{\al\bet\ga}x_i^\bet
x_j^\ga\n^\al\\
&-\dot\F F^{ij}\riema \al\bet\ga\de\n^\al x_i^\bet x_k^\ga x_j^\de\h_\e x_l^\e g^{kl}\\
&-\tilde f_\bet x_i^\bet x_k^\al \h_\al g^{ik}- \tilde f_{\n^\bet}x^\bet_k h^{ik}x^\al_i\h_\al,
\end{aligned}
\end{equation}
where $\h$ is the covariant vector field $(\h_\al)=e^{\psi}(-1,0,\dotsc,0)$.
\el
The proof is identical to the proof of \rl{1.4.4}.

In the Riemannian case we have:

\bl[Evolution of $v$]\lal{1.4.10}
Consider the flow \re{1.3.5} in a normal Gaussian coordinate system $(x^\al)$, where the $M(t)$ can be written as graphs of a function $u(t)$ over some compact Riemannian manifold $\so$. Then the quantity
\begin{equation}
v=\sqrt{1+\abs {Du}^2}=(r_\al\nu^\al)^{-1}
\end{equation}
satisfies the evolution equation
\begin{equation}\lae{1.4.4.10}
\begin{aligned}
\dot{ v}-\dot\F F^{ij} v_{ij}=&-\dot\F F^{ij}h_{ik}h_j^k v-2 v^{-1} \dot\F
F^{ij}v_iv_j \\
&+[(\F- f)-\dot\F F]r_{\al\bet}\n^\al\n^\bet v^2\\
&+2\dot\F F^{ij}h_j^k x_i^\al x_k^\bet r_{\al\bet} v^2+\dot\F F^{ij}r_{\al\bet\ga}x_i^\bet
x_j^\ga\n^\al v^2\\
&+\dot\F F^{ij}\riema \al\bet\ga\de\n^\al x_i^\bet x_k^\ga x_j^\de r_\e x_l^\e g^{kl} v^2\\
&+\tilde f_\bet x_i^\bet x_k^\al r_\al g^{ik}v^2+\tilde f_{\n^\bet}x^\bet_k h^{ik}x^\al_i
r_\al v^2,
\end{aligned}
\end{equation}
where $r=x^0$ and $(r_\al)=(1,0,\dotsc,0)$.
\el

\section{Existence results}

From now on we shall assume that ambient manifold $N$ is  Lorentzian, or more precisely, that it is smooth, globally hyperbolic with a  compact, connected Cauchy hypersurface. Then there exists a smooth future oriented time function $x^0$ such that the metric in $N$ can be expressed in Gaussian coordinates $(x^\al)$ as 
\begin{equation}\lae{6.1.1}
d\bar s^2=e^{2\psi}\{-(dx^0)^2+\s_{ij}dx^idx^j\},
\end{equation}
where $x^0$ is the time function and the $(x^i)$ are local coordinates for 
\begin{equation}
\so=\{x^0=0\}.
\end{equation}
$\so$ is then also a compact, connected Cauchy hypersurface. For a proof of the splitting result see \cite[Theorem 1.1]{sanchez:splitting}, and for the fact that all Cauchy hypersurfaces are diffeomorphic and hence $\so$ is also compact and connected, see \cite[Lemma 2.2]{Bernal:2003jb}.

One advantage of working in globally hyperbolic spacetimes with a compact Cauchy hypersurface is that all compact, connected  spacelike $C^m$-hypersurfaces $M$ can be written as graphs over $\so$.

\bl\lal{6.1.1}
Let $N$ be as above and $M\su N$ a connected, spacelike hypersurface of class $C^m$, $1\le m$, then $M$ can be written as a graph over $\so$
\begin{equation}
M=\graph \fv u{\so}
\end{equation}
with $u\in C^m(\so)$.
\el

We proved this lemma under the additional hypothesis that $M$ is achronal, \cite[Proposition 2.5]{cg:indiana},  however, this assumption is unnecessary as has been shown in \cite[Theorem 1.1]{heiko:diplom}.

We are looking at the curvature flow \fre{1.3.5} and want to prove that it converges to a stationary solution hypersurface, if certain assumptions are satisfied.

The existence proof consists of four steps:
\bir
\item
Existence on a maximal time interval $[0,T^*)$.

\item
Proof that the flow stays in a compact subset.

\item
Uniform a priori estimates in an appropriate function space, e.g., $C^{4,\al}(\so)$ or $C^{\un}(\so)$, which, together with (ii), would imply $T^*=\un$.

\item
Conclusion that the flow---or at least a subsequence of the flow hypersur\-faces---converges if $t$ tends to infinity.
\ei

The existence on a maximal time interval is always guaranteed, if the data are sufficiently regular, since the problem is parabolic. If the flow hypersurfaces can be written as graphs in a Gaussian coordinate system, as will always be the case in a globally hyperbolic spacetime with a compact Cauchy hypersurface in view of \rl{6.1.1}, the conditions are better than in the general case:
\bt\lat{6.2.2}
Let $4\le m\in\N$ and $0<\al<1$, and assume the semi-Riemannian space $N$ to be of class $C^{m+2,\al}$. Let the strictly monotone curvature function $F$, the functions $f$ and $\F$ be of class $C^{m,\al}$ and let $M_0\in C^{m+2,\al}$ be an admissible compact, spacelike, connected, orientable\footnote{Recall that oriented simply means there exists a continuous normal, which will always be the case in a globally hyperbolic spacetime.} hypersurface. Then the curvature flow \fre{1.3.5} with initial hypersurface $M_0$ exists in  a maximal time interval $[0,T^*)$, $0<T^*\le\un$, where in case that the flow hypersurfaces cannot be expressed as graphs they are supposed to be smooth, i.e, the conditions should be valid for arbitrary $4\le m\in\N$ in this case.
\et

A proof can be found in \cite[Theorem 2.5.19, Lemma 2.6.1]{cg:cp}.

\cvm
The second step, that the flow stays in a compact set, can only be achieved by barrier assumptions, \cf \rd{1.6.7}. Thus, let $\Om\su N$ be open and precompact such that $\pa\Om$ has exactly two components
\begin{equation}\lae{6.4.1}
\pa\Om=M_1\uud M_2
\end{equation}
where $M_1$ is a lower barrier for the pair $(F,f)$ and $M_2$ an upper barrier. Moreover, $M_1$ has to lie in the past of $M_2$
\begin{equation}\lae{6.5.1}
M_1\su I^-(M_2),
\end{equation}
\cf \cite[Remark 2.7.8]{cg:cp}.

Then the flow hypersurfaces will always stay inside $\bar\Om$, if the initial hypersurface $M_0$ satisfies $M_0\su \Om$, \cite[Theorem 2.7.9]{cg:cp}. This result is also valid if $M_0$ coincides with one the barriers, since then the velocity $(\F-\tilde f)$ has a weak sign and the flow moves into $\Om$ for small $t$, if it moves at all, and the arguments of the proof are applicable.

In Lorentzian manifolds the existence of barriers is associated with the presence of past and future singularities. In globally hyperbolic spacetimes, when $N$ is topologically a product
\begin{equation}
N=I\times\so,
\end{equation}
where $I=(a,b)$, singularities can only occur, when the endpoints of the interval are approached. A singularity, if one exists, is called a \tit{crushing singularity}, if the sectional curvatures become unbounded, i.e.,
\begin{equation}
\riema\al\bet\ga\de \bar R^{\al\bet\ga\de}\ra\un
\end{equation}
and such a singularity should provide a future \resp past barrier for the mean curvature function $H$.

\bd\lad{3.6.1.2}
Let $N$ be a globally hyperbolic spacetime with compact Cauchy hypersurface $\so$
so that $N$ can be written as a topological product $N=I\times \so$ and its
metric expressed as
\begin{equation}\lae{3.6.0.4}
d\bar s^2=e^{2\psi}(-(dx^0)^2+\s_{ij}(x^0,x)dx^idx^j).
\end{equation}
Here, $x^0$ is a globally defined future directed time function and $(x^i)$ are local
coordinates for $\so$.
$N$ is said to have a \tit{future} \resp \tit{past mean
curvature barrier}, if there are sequences $M_k^+$ \resp $M_k^-$ of closed,
spacelike, admissible  hypersurfaces such that
\begin{equation}\lae{6.9.1}
\lim_{k\ra\un} \fv H{M_k^+}=\un \q\tup{\resp}\q \lim_{k\ra\un} \fv H{M_k^-}=-\un
\end{equation}
and
\begin{equation}\lae{6.10.1}
\limsup \inf_{M_k^+}x^0> x^0(p)\qq\A\,p\in N
\end{equation}
\resp
\begin{equation}\lae{6.11.1}
\liminf \sup_{M_k^-}x^0< x^0(p)\qq\A\,p\in N,
\end{equation}
\ed

If one stipulates that the principal curvatures of the $M_k^+$ \resp $M_k^-$ tend to plus \resp minus infinity, then these hypersurfaces could also serve as barriers for other curvature functions. The past barriers would most certainly be non-admissible for any curvature function except $H$. 

\br\lar{6.4.2}
Notice that the assumptions \re{6.9.1} alone already implies \re{6.10.1} \resp \re{6.11.1}, if either
\begin{equation}\lae{6.12.1}
\limsup \inf_{M_k^+}x^0>a
\end{equation}
\resp
\begin{equation}
\liminf \sup_{M_k^-}x^0 <b
\end{equation}
where $(a,b)=x^0(N)$, or, if 
\begin{equation}\lae{6.14.1}
\bar R_{\al\bet}\nu^\al\nu^\bet\ge -\Lam\qq\A\, \spd\nu\nu=-1.
\end{equation}
where $\Lam\ge 0$. 
\er

\bp
It suffices to prove that the relation \re{6.10.1} is automatically satisfied under the assumptions \re{6.12.1} or \re{6.14.1}  by switching the light cone and replacing $x^0$ by $-x^0$ in case of the past barrier. 

Fix $k$, and let 
\begin{equation}
\tau_k=\inf_{M_k}x^0,
\end{equation}
 then the coordinate slice 
 \begin{equation}
M_{\tau_k}=\{x^0=\tau_k\}
\end{equation}
touches $M_k$ from below in a point $p_k\in M_k$ where $\tau_k=x^0(p_k)$ and the maximum principle yields that in that point
\begin{equation}
\fv H{M_{\tau_k}}\ge \fv H{M_k},
\end{equation}
hence, if $k$ tends to infinity the points $(p_k)$ cannot stay in a compact subset, i.e., 
\begin{equation}\lae{6.18.2}
\limsup x^0(p_k)\ra b
\end{equation}
or
\begin{equation}\lae{6.19.1}
\limsup x^0(p_k)\ra a.
\end{equation}

We shall show that only \re{6.18.2} can be valid. The relation \re{6.19.1} evidently contradicts \re{6.12.1}.

In case the assumption \re{6.14.1}  is valid, we consider a fixed coordinate slice $M_0=\{x^0=\const\}$, then all hypersurfaces $M_k$ satisfying
\begin{equation}
\fv H{M_0}<\inf_{M_k}H\q\wed\q \sqrt{n\Lam}<\inf_{M_k}H
\end{equation}
have to lie in the future of $M_0$, \cf \cite[Lemma 4.7.1]{cg:cp}, hence the result.
\ep

A future mean curvature barrier certainly represents a singularity, at least if $N$ 
satisfies the condition
\begin{equation}
\bar R_{\al\bet}\nu^\al\nu^\bet\ge -\Lam \qq\A\,\spd\nu\nu=-1
\end{equation}
where $\Lam\ge 0$, because of the future timelike incompleteness, which is proved in \cite{galloway:cft}, and is a generalization of Hawking's earlier result for $\Lam=0$, \cite{he:book}.  But these
singularities need not be crushing, \cf \cite[Section 2]{cg:imcf} for a counterexample.

\cvm
The uniform a priori estimates for the flow hypersurfaces are the hardest part in any existence proof. When the flow hypersurfaces can be written as graphs it suffices to prove $C^1$ and $C^2$ estimates, namely, the induced metric
\begin{equation}
g_{ij}(t,\xi)=\spd{x_i}{x_j}
\end{equation}
where $x=x(t,\xi)$ is a local embedding of the flow, should stay uniformly positive definite, i.e., there should exist positive constants $c_i$, $1\le i\le 2$, such that
\begin{equation}
c_1 g_{ij}(0,\xi)\le g_{ij}(t,\xi)\le c_2 g_{ij}(0,\xi),
\end{equation}
or equivalently, that the quantity
\begin{equation}
\tilde v=\spd{\h}\nu,
\end{equation}
where $\nu$ is the past directed normal of $M(t)$ and $\h$ the vector field
\begin{equation}
\h=(\h_\al)=e^\psi(-1,0,\ldots,0),
\end{equation}
is uniformly bounded, which is achieved with the help of the parabolic equation \fre{1.4.4.3}, if it is possible at all.

However, in some special situations $C^1$-estimates are automatically satisfied, \cf \rt{1.6.11} at the end of this section.

For the $C^2$-estimates  the principle curvatures $\ka_i$ of the flow hypersurfaces have to stay in a compact set in the cone of definition $\C$ of $F$, e.g., if $F$ is the Gaussian curvature, then $\C=\C_+$ and one has to prove that there are positive constants $k_i$, $i=1,2$ such that
\begin{equation}
k_1\le \ka_i\le k_2\qq\A\, 1\le i\le n
\end{equation}
uniformly in the cylinder $[0,T^*)\times M_0$, where $M_0$ is  any manifold that can serve as a base manifold for the embedding $x=x(t,\xi)$.

The parabolic equations that are used for these curvature estimates are \fre{1.4.3.13}, usually for an upper estimate, and \fre{1.4.3.10} for the lower estimate. Indeed, suppose that the flow starts at the upper barrier, then
\begin{equation}\lae{6.18.1}
F\ge f
\end{equation}
at $t=0$ and this estimate remains valid throughout the evolution because of the parabolic maximum principle, use \re{1.4.3.13}. Then, if upper estimates for the $\ka_i$ have been derived and if $f>0$ uniformly, then we conclude from \re{6.18.1} that the $\ka_i$ stay in a compact set inside the open cone $\C$, since 
\begin{equation}
\fv F{\pa\C}=0.
\end{equation}

To obtain higher order estimates we are going to exploit the fact that the flow hypersurfaces are graphs over $\so$ in an essential way, namely, we look at the associated scalar flow equation \fre{1.4.21} satisfied by $u$. This equation is a nonlinear uniformly parabolic equation, where the operator $\F(F)$ is also concave in $h_{ij}$, or equivalently, convex in $u_{ij}$, i.e., the $C^{2,\al}$-estimates of Krylov and Safonov, \cite[Chapter 5.5]{nk} or see \cite[Chapter 10.6]{oliver:pde} for a very clear and readable presentation, are applicable, yielding uniform estimates for the standard parabolic H\"older semi-norm
\begin{equation}\lae{6.29.2}
[D^2 u]_{\bet,\bar Q_T}
\end{equation}
for some $0<\bet\le\al$ in the cylinder
\begin{equation}
Q=[0,T)\times \so,
\end{equation}
independent of $0<T<T^*$, which in turn will lead to $H^{m+2+\al,\frac{m+2+\al}2}(\bar Q_T)$ estimates, \cf \cite[Theorem 2.5.9, Remark 2.6.2]{cg:cp}.  

$H^{m+2+\al,\frac{m+2+\al}2}(\bar Q_T)$ is a parabolic H\"older space, \cf \cite[p. 7]{lus} for the original definition and \cite[Note 2.5.4]{cg:cp} in the present context.

The estimate \re{6.29.2} combined with the uniform $C^2$-norm leads to uniform
 $C^{2,\bet}(\so)$-estimates independent of $T$.

These estimates imply that $T^*=\un$.

\cvm
Thus, it remains to prove that $u(t,\cdot)$ converges in $C^{m+2}(\so)$ to a stationary solution $\tilde u$, which is then also of class $C^{m+2,\al}(\so)$ in view of the Schauder theory.

Because of the preceding a priori estimates $u(t,\cdot)$ is precompact in $C^2(\so)$. Moreover, we deduce from the scalar flow equation \fre{1.4.21} that $\dot u$ has a sign, i.e., the $u(t,\cdot)$ converge monotonely in $C^0(\so)$ to $\tilde u$ and therefore also in $C^2(\so)$.

To prove that $\graph \tilde u$ is a solution, we again look at \re{1.4.21} and integrate it with respect to $t$ to obtain for fixed $x\in\so$
\begin{equation}
\abs{\tilde u(x)-u(t,x)}=\int_t^\un e^{-\psi}v\abs{\F-\tilde f},
\end{equation}
where we used that $(\F-\tilde f)$ has a sign, hence $(\F-\tilde f)(t,x)$ has to vanish when $t$ tends to infinity, at least for a subsequence, but this suffices to conclude that $\graph \tilde u$ is a stationary solution and
\begin{equation}
\lim_{t\ra\un}(\F-\tilde f)=0.
\end{equation}

Using the convergence of $u$ to $\tilde u$ in $C^2$, we can then prove:
\bt\lat{6.5.5}
The functions $u(t,\cdot)$ converge in $C^{m+2}(\so)$ to $\tilde u$, if the data satisfy the assumptions in \rt{6.2.2}, since we have
\begin{equation}\lae{6.33.2}
u\in H^{m+2+\bet,\frac{m+2+\bet}2}(\bar Q),
\end{equation}
where $Q=Q_\un$. 
\et

\bp
Out of convention let us write $\al$ instead of $\bet$ knowing that $\al$ is the H\"older exponent in \re{6.29.2}.

We shall reduce the Schauder estimates to the standard Schauder estimates in $\R[n]$ for the heat equation with a right-hand side by using the already established results \re{6.29.2} and 
\begin{equation}
u(t,\cdot)\underset{C^2(\so)}\ra \tilde u\in C^{m+2,\al}(\so).
\end{equation}

Let $(U_k)$ be a finite open covering of $\so$ such that each $U_k$ is contained in a coordinate chart and
\begin{equation}\lae{6.35.2}
\diam U_k<\rho,
\end{equation}
$\rho$ small, $\rho$ will be specified in the proof, and let $(\h_k)$ be a subordinate finite partition of unity of class $C^{m+2,\al}$.

Since
\begin{equation}
u\in H^{m+2+\al,\frac{m+2+\al}2}(\bar Q_T)
\end{equation}
for any finite $T$, \cf \cite[Lemma 2.6.1]{cg:cp}, and hence 
\begin{equation}
u(t,\cdot)\in C^{m+2,\al}(\so)\qq\A\,0\le t<\un
\end{equation}
we shall choose $u_0=u(t_0,\cdot)$ as initial value for some large $t_0$ such that
\begin{equation}\lae{6.38.2}
\abs{a^{ij}(t,\cdot)-\tilde a^{ij}}_{0,\so}<\e/2\qq\A\,t\ge t_0,
\end{equation}
where
\begin{equation}
a^{ij}=v^2\dot\F F^{ij}
\end{equation}
and  $\tilde a^{ij}$ is defined correspondingly for $\tilde M=\graph \tilde u$.

However, making a variable transformation we shall always assume that $t_0=0$ and $u_0=u(0,\cdot)$.

We shall prove \re{6.33.2} successively.

\cvm
(i) Let us first show that
\begin{equation}
D_xu\in H^{2+\al,\frac{2+\al}2}(\bar Q).
\end{equation}
This will be achieved, if we show that for an arbitrary $\xi\in C^{m+1,\al}(T^{1,0}(\so))$
\begin{equation}
\f=D_\xi u\in H^{2+\al,\frac{2+\al}2}(\bar Q),
\end{equation}
\cf \cite[Remark 2.5.11]{cg:cp}.

Differentiating the scalar flow equation \fre{1.4.21} with respect to $\xi$ we obtain
\begin{equation}\lae{6.42.2}
\dot \f-a^{ij}\f_{ij}+b^i\f_i+c\f=f,
\end{equation}
where of course the symbol $f$ has a different meaning then in \re{1.4.21}.

Later we want to apply the Schauder estimates for solutions of the heat flow equation with right-hand side. In order to use elementary potential estimates we have to cut off $\f$ near the origin $t=0$ by considering
\begin{equation}
\tilde\f=\f \theta,
\end{equation}
where $\theta=\theta(t)$ is smooth satisfying
\begin{equation}
\theta(t)=
\begin{cases}
1,&t>1,\\
0,&t\le\frac12.
\end{cases}
\end{equation}
This modification doesn't cause any problems, since we already have a priori estimates for finite $t$, and we are only concerned about the range $1\le t<\un$. $\tilde\f$ satisfies the same equation as $\f$ only the right-hand side has the additional summand $w\dot\theta$.

Let $\h=\h_k$ be one of the members of the partition of unity and set
\begin{equation}
w=\tilde\f\h,
\end{equation}
then $w$ satisfies a similar equation with slightly different right-hand side
\begin{equation}\lae{6.44.2}
\dot w-a^{ij}w_{ij}+b^iw_i+cw=\tilde f
\end{equation}
but we shall have this in mind when applying the estimates.

The $w(t,\cdot)$ have compact support in one of the $U_k$'s, hence we can replace the covariant derivatives of $w$ by ordinary partial derivatives without changing the structure of the equation and the properties of the right-hand side, which still only depends linearly on $\f$ and $D\f$.

We want to apply the well-known estimates for the ordinary heat flow equation
\begin{equation}
\begin{aligned}
\dot w-\D w&=\hat f
\end{aligned}
\end{equation}
where $w$ is defined  in $\R[]\times\R[n]$.

To reduce the problem to this special form, we pick an arbitrary $x_0\in U_k$, set $z_0=(0,x_0)$, $z=(t,x)$ and consider instead of \re{6.44.2}
\begin{equation}
\begin{aligned}
\dot w-a^{ij}(z_0)w_{ij}&=\hat f\\
&=[a^{ij}(z)-a^{ij}(z_0)]w_{ij}-b^iw_i-cw+\tilde f,
\end{aligned}
\end{equation}
where we emphasize that the difference
\begin{equation}\lae{6.47.2}
\abs{a^{ij}(z)-a^{ij}(z_0)}
\end{equation}
can be made smaller than any given $\e>0$ by choosing $\rho=\rho(\e)$ in \re{6.35.2} and $t_0=t_0(\e)$ in \re{6.38.2} accordingly. Notice also that this equation can be extended into $\R[]\times \R[n]$, since all functions have support in $\{t\ge \frac12\}$.

Let $0<T<\un$ be arbitrary, then all terms belong to the required function spaces in $\bar Q_T$ and there holds
\begin{equation}
[w]_{2+\al,Q_T}\le c[\hat f]_{\al,Q_T},
\end{equation}
where $c=c(n,\al)$. The brackets indicate the standard unweighted parabolic semi-norms, \cf \cite[Definition 2.5.2]{cg:cp}, which are identical to those defined in \cite[p. 7]{lus}, but there the brackets are replaced by kets. 

Thus, we conclude
\begin{equation}
\begin{aligned}
\hp{a}&[w]_{2+\al,Q_T}\le c \sup_{U_k\times (0,T)}\abs{a^{ij}(z)-a^{ij}(z_0)}[D^2w]_{\al,Q_T}+c[f]_{\al,Q_T} \\
&+c_1\{[D^2u]_{\al,Q_+}+[Du]_{\al,Q_T}+[u]_{\al,Q_T}+\abs{w}_{0,Q_T}+\abs{D^2w}_{0,Q_T}\},
\end{aligned}
\end{equation}
where $c_1$ is independent of $T$, but dependent on $\h_k$. Here we also used the fact that the lower order coefficients and $\f, D\f$ are uniformly bounded.

Choosing now $\e>0$ so small that
\begin{equation}
c\msp \e<\tfrac12
\end{equation}
and $\rho$, $t_0$ accordingly such the difference in \re{6.47.2} is smaller than $\e$, we deduce
\begin{equation}
\begin{aligned}
\hp{a}&[w]_{2+\al,Q_T}\le 2c[f]_{\al,Q_T} \\
& +2c_1\{[D^2u]_{\al,Q_+}+[Du]_{\al,Q_T}+[u]_{\al,Q_T}+\abs{w}_{0,Q_T}+\abs{D^2w}_{0,Q_T}\}.
\end{aligned}
\end{equation}

Summing over the partition of unity and noting that $\xi$ is arbitrary we see that in the preceding inequality we can replace $w$ by $Du$ everywhere  resulting in the estimate
\begin{equation}
\begin{aligned}
&\hp{a}[Du]_{2+\al,Q_T}\le 
c_1[f]_{\al,Q_T}\\
&\qq +c_1\{[D^2u]_{\al,Q_T}+[Du]_{\al,Q_T}+[u]_{\al,Q_T}+\abs{Du}_{0,Q_T}+\abs{D^3u}_{0,Q_T}\},
\end{aligned}
\end{equation}where $c_1$ is a new constant still independent of $T$.

Now the only critical terms on the right-hand side are $\abs{D^3u}_{0,Q_T}$, which can be estimated by \re{6.55.2}, and the H\"older semi-norms with respect to $t$
\begin{equation}
[Du]_{\frac\al 2,t,Q_T}+[u]_{\frac\al 2,t,Q_T}.
\end{equation}
The second  one is taken care of by the boundedness of $\dot u$, see \fre{1.4.21}, while the first one is estimated with the help of equation \re{6.42.2} revealing
\begin{equation}
\abs{D\dot u}\le c\{\sup_{[0,T]}\abs{u}_{3,\so}+\abs{f}_{0,Q_T}\},
\end{equation}
since for fixed but arbitrary $t$ we have
\begin{equation}\lae{6.55.2}
\abs{u}_{3,\so}\le \e[D^3u]_{\al,\so}+c_\e\abs{u}_{0,\so},
\end{equation}
where $c_\e$ is independent of $t$.

Hence we conclude
\begin{equation}
\abs{Du}_{2+\al,Q_T}\le\const
\end{equation}
uniformly in $T$.

\cvm
(ii) Repeating these estimates successively for $2\le l\le m$ we obtain uniform estimates for
\begin{equation}
\sum_{l=2}^m[D^l_xu]_{2+\al,Q_T},
\end{equation}
which, when combined with the uniform $C^2$-estimates, yields
\begin{equation}\lae{6.58.2}
\abs{u(t,\cdot)}_{m+2,\al,\so}\le \const
\end{equation}
uniformly in $0\le t<\un$.

Looking at the equation \re{1.4.21} we then deduce
\begin{equation}\lae{6.59.2}
\abs{\dot u(t,\cdot)}_{m,\al,\so}\le \const
\end{equation}
uniformly in $t$.

\cvm
(iii) To obtain the estimates for $D^r_tu$ up to the order
\begin{equation}
[\tfrac{m+2+\al}2]
\end{equation}
we differentiate the scalar curvature  equation with respect to $t$ as often as necessary and also with respect to the mixed derivatives $D^r_tD^s_x$ to estimate
\begin{equation}
\sum_{1\le 2r+s<m+2+\al}D^r_tD^s_x u
\end{equation}
using \re{6.58.2}, \re{6.59.2} and the results from the prior differentiations. 

Combined with the estimates for the heat equation in $\R[]\times\R[n]$ these estimates will also yield the necessary a priori estimates for the H\"older semi-norms in $\bar Q$, where again the smallness of \re{6.47.2} has to be used repeatedly.
\ep

\br
The preceding regularity result is also valid in Riemannian manifolds, if the flow hypersurfaces can be written as graphs in a Gaussian coordinate system. In fact the proof is unaware of the nature of the ambient space.
\er

With the method described above the following existence results have been proved in globally hyperbolic spacetimes with a compact Cauchy hypersurface. $\Om\su N$ is always a precompact domain the boundary of which is decomposed as in \re{6.4.1} and \re{6.5.1} into an upper and lower barrier for the pair $(F,f)$. We also apply the stability results from \cite[Section 5]{cg:survey} and the just proved regularity of the convergence and formulate the theorems accordingly. 

By convergence of the flow in $C^{m+2}$ we mean convergence of the leaves $M(t)=\graph u(t,\cdot)$ in this norm.

\bt\lat{5.7.0.0}
Let $M_1$, $M_2$ be lower \resp upper barriers for the pair $(H,f)$, where $f\in C^{m,\al}(\bar\Om)$ and the $M_i$ are of class $C^{m+2,\al}$, $4\le m$, $0<\al<1$, then the curvature flow
\begin{equation}
\begin{aligned}
\dot x&=(H-f)\nu\\
x(0)&=x_0,
\end{aligned}
\end{equation}
where $x_0$ is an embedding of the initial hypersurface $M_0=M_2$ exists for all time and converges in $C^{m+2}$ to a stable solution $M$ of class $C^{m+2,\al}$ of the equation
\begin{equation}
\fv HM=f,
\end{equation}
provided the initial hypersurface is not already a solution.
\et

The existence result was proved in \cite[Theorem 2.2]{cg:mz}, see also \cite[Theorem 4.2.1]{cg:cp} and the remarks following the theorem, and the stability result in \cite[Section 5]{cg:survey}. Notice that $f$ isn't supposed to satisfy any sign condition.

For spacetimes that satisfy the timelike convergence condition and for functions $f$ with special structural conditions existence results via a mean curvature flow were first proved in \cite{eh1}.

The Gaussian curvature or the curvature functions $F$ belonging to the larger class $(K^*)$, see \cite{cg:indiana} for a definition, require that the admissible hypersurfaces are strictly convex. 

Moreover, proving a priori estimates for the second fundamental form of a hypersurface $M$ in general semi-Riemannian manifolds, when the curvature function is not the mean curvature, or does not behave similar to it,  requires that a strictly convex function $\chi$ is defined in a neighbourhood of the hypersurface, see \frl{2.8.1} where sufficient assumptions are stated which imply the existence of strictly convex  functions.

Furthermore, when we consider curvature functions of class $(K^*)$, notice that the Gaussian curvature belongs to that class, then the right-hand side $f$ can be defined in $T(\bar\Om)$ instead of $\bar\Om$, i.e., in a local trivialization of the tangent bundle $f$ can be expressed as
\begin{equation}
f=f(x,\nu)\q\wed\q\nu\in T_x(N).
\end{equation}
We shall formulate the existence results with this more general assumption.

\bt
Let $F\in C^{m,\al}(\C_+)$, $4\le m$, $0<\al<1$, be a curvature function of class $(K^*)$, let $0<f\in C^{m,\al}(T(\bar\Om))$, and let $M_1$, $M_2$ be lower \resp upper barriers for $(F,f)$ of class $C^{m+2,\al}$. Then the curvature flow
\begin{equation}
\begin{aligned}
\dot x&=(\F-\tilde f)\nu\\
x(0)&=x_0
\end{aligned}
\end{equation}
where $\F(r)=\log r$ and $x_0$ is an embedding of $M_0=M_2$, exists for all time and converges in $C^{m+2}$ to a stationary solution $M\in C^{m+2,\al}$ of the equation
\begin{equation}
\fv FM=f
\end{equation}
provided the initial hypersurface $M_2$ is not already a stationary solution and there exists a strictly convex function $\chi\in C^2(\bar\Om)$.
\et

The theorem was proved in \cite{cg:indiana} when $f$ is only defined in $\bar\Om$ and in the general case in \cite[Theorem 4.1.1]{cg:cp}.

\cvm
When $F=H_2$ is the scalar curvature operator, then the requirement that $f$ is defined in the tangent bundle and not merely in $N$ is a necessity, if the scalar curvature is to be prescribed. To prove existence results in this case, $f$ has to satisfy some natural structural conditions, namely,
\begin{align}
0<c_1&\le f(x,\n)\qq\tup{if}\q\spd\n\n=-1,\lae{0.3}\\[\cma]
\nnorm{f_\bet(x,\n)}&\le c_2 (1+\nnorm\n^2),\lae{0.4}\\
\intertext{and}
\nnorm{f_{\n^\bet}(x,\n)}&\le c_3 (1+\nnorm\n),\lae{0.5}
\end{align}
\nd for all $x\in\bar\Om$ and all past directed timelike vectors $\n\in T_x(\Om)$,
where $\nnorm{\cdot}$ is a Riemannian reference metric.

Applying a curvature flow to obtain stationary solutions requires to approximate  $f$ by functions  $f_k$ and to use these functions for the flow.  

The functions $f_k$ have the property that $\nnorm{f_{k\bet}}$ only grows linearly in $\nnorm\nu$ and $\nnorm{f_{k\n^\bet}(x,\n)}$ is bounded. To simplify the presentation we shall therefore assume that $f$ satisfies
\begin{equation}\lae{6.33.1}
\nnorm{f_\bet(x,\n)}\le c_2 (1+\nnorm\n),
\end{equation}
\begin{equation}\lae{6.34.1}
\nnorm{f_{\n^\bet}(x,\n)}\le c_3,
\end{equation}
and also 
\begin{equation}\lae{6.35.1}
0<c_1\le f(x,\n)\qq\A\,\nu\in T_x(N),\;\spd\nu\nu<0,
\end{equation}
although the last assumption is only a minor point that can easily be dealt with, see \cite[Remark 2.6]{cg:scalar}, and \cite[Section 7 and 8]{cg:scalar} for the other approximations of $f$.

Now, we can formulate the existence result for the scalar curvature operator $F=H_2$ under these provisions.
\bt\lat{6.6.1}
Let $f\in C^{m,\al}(T(\bar\Om))$, $4\le m$, $0<\al<1$, satisfy the conditions \re{6.33.1}, \re{6.34.1} and \re{6.35.1}, and let $M_1$, $M_2$ be lower \resp upper barriers of class $C^{m+2,\al}$ for $(F,f)$. Then the curvature flow for $F$
\begin{equation}
\begin{aligned}
\dot x&=(\F-\tilde f)\\
x(0)&=x_0
\end{aligned}
\end{equation}
where $\F(r)=r^\frac12$ and $x_0$ is an embedding of $M_0=M_2$, exists for all time and converges in $C^{m+2}$ to a stationary solution $M\in C^{m+2,\al}$ of
\begin{equation}\lae{c5.76}
\fv{F}{M}=f
\end{equation}
provided there exists a strictly convex function $\chi\in C^2(\bar\Om)$.
\et

This theorem has been proved by Christian Enz in \cite{enz:scalar} using the curvature estimates in \cite{cg:weingarten}, see also \frs{h6}.

The first existence result for equation \re{c5.76} was proved in \cite{cg:scalar} by considering $\e$-regularizations of the scalar curvature function in the curvature flow and by proving rather elaborate curvature estimates. The new existence proof is much simpler and more elegant.

To conclude this section let us show which spacelike hypersurfaces satisfy $C^1$-estimates automatically.

\bt\lat{1.6.11}
Let $M=\graph \fv u\so$ be a compact, spacelike hypersurface represented in a
Gaussian coordinate system with  unilateral bounded  principal curvatures,
e.g.,
\begin{equation}
\kappa_i\ge \kappa_0\q\A\,i.
\end{equation}
Then, the quantity $\tilde v=\frac{1}{\sql}$ can be estimated by
\begin{equation}
\tilde v\le c(\abs u,\so,\sigma_{ij},\psi,\kappa_0),
\end{equation}
where we assumed that in the Gaussian coordinate system the
ambient metric has the form as in \re{6.1.1}.
\et
\bp
We suppose as usual that the Gaussian coordinate system is future oriented, and that the
second fundamental form is evaluated with respect to the past directed normal.
We observe that
\begin{equation}\lae{1.6.87}
\norm{Du}^2=g^{ij}u_iu_j=e^{-2\psi}\frac{\abs{Du}^2}{v^2}\raise 2pt \hbox{,}
\end{equation}
hence, it is equivalent to find an a priori estimate for $\norm{Du}$.

Let $\lambda$ be a real parameter to be specified later, and set
\begin{equation}
w=\tfrac{1}{2}\log\norm{Du}^2+\lambda u.
\end{equation}
We may regard $w$ as being defined on $\so$; thus, there is $x_0\in\so$ such that
\begin{equation}
w(x_0)=\sup_\so w,
\end{equation}
and we conclude
\begin{equation}
0=w_i=\frac{1}{\norm{Du}^2}\,u_{ij}u^j+\lambda u_i
\end{equation}
in $x_0$, where the covariant derivatives are taken with respect to the induced
metric
$g_{ij}$, and the indices are also raised with respect to that metric.

Expressing the second fundamental form of a graph with the help of the Hessian of the function
\begin{equation}\lae{2.16}
e^{-\psi}v^{-1}h_{ij}=-u_{ij}-\cha 000\mspace{1mu}u_iu_j-\cha 0i0
\mspace{1mu}u_j-\cha 0j0\mspace{1mu}u_i-\cha ij0.
\end{equation}
 we deduce further
\begin{equation}\lae{4.16}
\begin{aligned}
\lambda\norm{Du}^4&=-u_{ij}u^iu^j\\
&= e^{-\psi}\tilde vh_{ij}u^iu^j+\cha 000\msp \norm{Du}^4\\
&\hp{=}\msp[2]+2\cha 0j0\msp u^j\norm{Du}^2+\cha ij0\msp u^iu^j.
\end{aligned}
\end{equation}
Now, there holds
\begin{equation}
u^i=g^{ij}u_j=e^{-2\psi}\sigma^{ij}u_jv^{-2},
\end{equation}
and by assumption,
\begin{equation}
h_{ij}u^iu^j\ge \kappa_0\msp\norm{Du}^2,
\end{equation}
i.e., the critical terms on the right-hand side of \re{4.16} are of fourth order in
$\norm{Du}$ with bounded coefficients, and we conclude that $\norm{Du}$ can't be
too large in $x_0$ if we choose $\lambda$ such that
\begin{equation}
\lambda\le -c\msp\nnorm{\cha \alpha\beta 0}-1
\end{equation}
with a suitable constant $c$; $w$, or equivalently, $\norm{Du}$ is therefore
uniformly bounded from above.
\ep

Especially for convex graphs over $\so$ the term $\tilde v$ is uniformly bounded as long as
they stay in a compact set.  

\section{Curvature flows in Riemannian manifolds}\las{h6}
The existence results for solutions of an equation like
\begin{equation}\lae{h1.1}
\fv FM=f
\end{equation}
or the long time existence of curvature flows, relies on a priori estimates. The derivation of second order estimates, i.e., curvature estimates, is usually the the most difficult task and depends crucially on the curvature functions involved. However, if the ambient space is Riemannian and the right-hand side, or external force, $f$ only depends on $x$, $f=f(x)$, then curvature estimates can be derived for general concave curvature functions provided lower order a priori estimates are already known, as we shall show in the following. These estimates were first proved  in \cite{cg:weingarten}.

When proving a priori estimates for solutions of \re{h1.1} the concavity of $F$ plays a central role. As usual we consider $F$ to be defined in a cone $\C$ as well as on the space of admissible tensors such that 
\begin{equation}
F(h_{ij})=F(\ka_i).
\end{equation}
Notice that curvature functions are always assumed to be symmetric and if $F\in C^{m,\al}(\C)$, $2\le m$, $0<\al<1$, then $F\in C^{m,\al}(\mc S_\C)$, where $\mc S_\C\su T^{0,2}(M)$ is the open set of admissible symmetric tensors with respect to the given metric $g_{ij}$. The result is due to Ball, \cite{ball}, see also \cite[Theorem 2.1.8]{cg:cp}.

The second derivatives of $F$ then satisfy
\begin{equation}
\begin{aligned}
F^{ij,kl}\h_{ij}\h_{kl}=\sum_{i,j}\frac{\pa^2F}{\pa\ka_i\pa\ka_j}\h_{ii}\h_{jj}+\sum_{i\ne j} \frac{F_i-F_j}{\ka_i-\ka_j}(\h_{ij})^2\le 0\;\A\, \h\in \mc S,
\end{aligned}
\end{equation}
where $\mc S\su T^{0,2}(M)$ is the space of symmetric tensors, if $F$ is concave in $\C$, \cf \cite[Lemma 1.1]{cg96}.

However, a mere non-positivity of the right-hand side is in general not sufficient to prove a priori estimates for the $\ka_i$ resulting in the fact that only for special curvature functions for which a stronger estimate was known such a priori estimates could be derived and the problem \re{h1.1} solved, if further assumptions are satisfied. 

Sheng et al.\ then realized in \cite{urbas:duke} that the term
\begin{equation}
\sum_{i\ne j} \frac{F_i-F_j}{\ka_i-\ka_j}(\h_{ij})^2
\end{equation}
was all that is needed to obtain the stronger concavity estimates under certain circumstances. Indeed, if the $\ka_i$ are labelled 
\begin{equation}\lae{h1.5}
\ka_1\le\cdots\le \ka_n,
\end{equation}
then there holds:
\bl
Let $F$ be concave and monotone, and assume $\ka_1<\ka_n$, then
\begin{equation}\lae{h1.8}
\sum_{i\ne j} \frac{F_i-F_j}{\ka_i-\ka_j}(\h_{ij})^2\le \frac2{\ka_n-\ka_1}\sum_{i=1}^n(F_n-F_i)(\h_{ni})^2
\end{equation}
for any symmetric tensor $(\h_{ij})$, where we used coordinates such that $g_{ij}=\de_{ij}$.
\el

\bp
Without loss of generality we may assume that the $\ka_i$ satisfy the strict inequalities
\begin{equation}
\ka_1<\cdots<\ka_n,
\end{equation}
since these points are dense. The concavity of $F$ implies
\begin{equation}
F_1\ge\cdots\ge F_n,
\end{equation}
\cf \cite[Lemma 2]{eh2}, where
\begin{equation}
F_i=\pde F{\ka_i}>0;
\end{equation}
the last inequality is the definition of monotonicity. 
The inequality then follows immediately.
\ep

The right-hand side of inequality \re{h1.8} is exactly the quantity that is needed to balance a bad technical term in the a priori estimate for $\ka_n$, at least in Riemannian manifolds, as we shall prove. Unfortunately, this doesn't work in Lorentzian spaces, because of a sign difference in the Gau{\ss} equations.

The assumptions on the curvature function are very simple.
\bas\laas{h1.2}
Let $\C\su\R[n]$ be an open, symmetric, convex cone containing $\C_+$ and let $F\in C^{m,\al}(\C)\ii C^0(\bar\C)$, $m\ge 4$, be  symmetric, monotone, homogeneous of degree $1$, and concave such that
\begin{equation}
F>0\qq\text{in}\q\C
\end{equation}
and 
\begin{equation}\lae{h1.11}
\fv F{\pa\C}=0.
\end{equation}
\eas

These conditions on the curvature function will suffice. They could have been modified, even relaxed, e.g.,  by  only requiring that $\log F$ is concave, but then the condition 
\begin{equation}
F^{ij}g_{ij}\ge c_0>0,
\end{equation}
which automatically holds, if $F$ is concave and homogeneous of degree  $1$, would have been added, destroying the aesthetic simplicity of \ras{h1.2}.

Our estimates apply equally well to solutions of an equation as well as to solutions of curvature flows. Since curvature flows encompass equations, let us state the main estimate for curvature flows.

Let $\Om\su N$ be precompact and connected, and $0<f\in C^{m,\al}(\bar\Om)$. We consider the curvature flow
\begin{equation}\lae{h1.13} 
\begin{aligned}
\dot x&=-(\F-\tilde f)\nu\\
x(0)&=x_0,
\end{aligned}
\end{equation}
where $\F$ is  $\F(r)=r$ and $\tilde f=f$, $x_0$ is the embedding of an initial admissible hypersurface $M_0$ of class $C^{m+2,\al}$ such that 
\begin{equation}
\F-\tilde f\ge 0\qq\text{at}\q t=0,
\end{equation}
where of course $\F=\F(F)=F$. We introduce the technical function $\F$ in the present case only to make a comparison with the formulas and results in the previous sections, which all use the notation for the more general flows, easier.

We assume that $\bar\Om$ is covered by a Gaussian coordinate system $(x^\al)$, $0\le 1\le n$, such that the metric can be expressed as 
\begin{equation}
d\bar s^2=e^{2\psi}\{(dx^0)^2+\s_{ij}dx^idx^j\}
\end{equation}
and $\bar\Om$ is covered by the image of the cylinder
\begin{equation}
I\times \so
\end{equation}
where $\so$ is a compact Riemannian manifold and $I=x^0(\bar\Om)$,  $x^0$ is a global coordinate defined in $\bar\Om$ and $(x^i)$ are local coordinates of $\so$.

Furthermore we assume that $M_0$ and the other flow hypersurfaces can be written as graphs over $\so$. The flow should exist in a maximal time interval $[0,T^*)$, stay in $\Om$, and uniform $C^1$-estimates should already have been established.

\br
The assumption on the existence of the Gaussian coordinate system and the fact that the hypersurfaces can be written as graphs could be replaced by assuming the existence of a unit vector field $\h\in C^2(T^{0,1}(\bar\Om))$ and of a constant $\theta>0$ such that
\begin{equation}\lae{h1.17}
\spd\h\nu\ge 2\theta
\end{equation}
uniformly during the flow, since this assumption would imply uniform $C^1$-estimates, which are the requirement that the induced metric can be estimated accordingly by  controlled metrics from below and above, and because the existence of such a vector field is essential for the curvature estimate.

If the flow hypersurfaces are graphs in a Gaussian coordinate system, then such a vector field is given by
\begin{equation}\lae{h1.18}
\h=(\h_\al)=e^{\psi}(1,0,\dots,0)
\end{equation}
and the $C^1$-estimates are tantamount to the validity of inequality \re{h1.17}.

In case $N=\R$ and starshaped hypersurfaces one could also use the term
\begin{equation}
\spd x\nu,
\end{equation}
\cf \cite[Lemma 3.5]{cg90}.
\er

Then we shall prove: 
\bt\lat{h1.4}
Under the assumptions stated above the principal curvatures $\ka_i$ of the flow hypersurfaces are uniformly bounded from above
\begin{equation}
\ka_i\le c,
\end{equation}
provided there exists a strictly convex function $\chi\in C^2(\bar\Om)$. The constant $c$ only depends on $\abs{f}_{2,\Om}$, $\theta$, $F(1,\dots,1)$, the initial data, and the estimates for $\chi$  and those of the ambient Riemann curvature tensor in $\bar\Om$.

Moreover, the $\ka_i$ will stay in a compact set of $\C$.
\et

As an application of this estimate our former results on the existence of a strictly convex hypersurface $M$ solving the equation \re{h1.1}, \cite{cg96,cg97}, which we proved for curvature functions $F$ of class $(K)$, are now valid for curvature functions $F$ satisfying \ras{h1.2} with $\C=\C_+$. 

We are even able to solve the existence problem by using a curvature flow which formerly only worked in case that the sectional curvature of the ambient space was non-positive.
\bt\lat{h1.5}
Let $F$ satisfy the assumptions above with $\C=\C_+$ and assume that the boundary of $\Om$ has two components
\begin{equation}
\pa\Om=M_1\uud M_2,
\end{equation}
where the $M_i$ are closed, connected strictly convex hypersurfaces of class $C^{m+2,\al}$, $m\ge 4$, which can be written as graphs in a normal Gaussian coordinate system covering $\bar\Om$, and where we assume that the normal of $M_1$ points outside of $\Om$ and that of $M_2$ inside. Let $0<f\in C^{m,\al}(\bar\Om)$,  and assume that $M_1$ is a lower barrier for the pair $(F,f)$ and $M_2$ an upper barrier, then the problem \re{h1.1} has a strictly convex solution $M\in C^{m+2,\al}$ provided there exists a strictly convex function $\chi\in C^2(\bar\Om)$. The solution is the limit hypersurface of a converging curvature flow.
\et

\subsection{Curvature estimates} 
Let $M(t)$ be the flow hypersurfaces, then their second fundamental form $h^j_i$ satisfies the evolution equation, \cf  \fre{1.4.1}: 
\bl\lal{h1.4.1}
The mixed tensor $h_i^j$ satisfies the parabolic equation 
\begin{equation}\lae{h1.4.1}
\begin{aligned}
&\qq\qq\dot h_i^j-\dot\F F^{kl}h_{i;kl}^j=\\[\cma]
&\hp{=}\; \dot\F F^{kl}h_{rk}h_l^rh_i^j-\dot\F F
h_{ri}h^{rj}+ (\F-\tilde f) h_i^kh_k^j\\
&\hp{+}-\tilde f_{\alpha\beta} x_i^\alpha x_k^\beta g^{kj}+ \tilde f_\alpha\n^\alpha h_i^j+\dot\F
F^{kl,rs}h_{kl;i}h_{rs;}^{\hphantom{rs;}j}\\
&\hp{=}+\ddot \F F_i F^j+2\dot \F F^{kl}\riema \alpha\beta\gamma\delta x_m^\alpha x_i ^\beta x_k^\gamma
x_r^\delta h_l^m g^{rj}\\
&\hp{=}-\dot\F F^{kl}\riema \alpha\beta\gamma\delta x_m^\alpha x_k ^\beta x_r^\gamma x_l^\delta
h_i^m g^{rj}-\dot\F F^{kl}\riema \alpha\beta\gamma\delta x_m^\alpha x_k ^\beta x_i^\gamma x_l^\delta h^{mj} \\
&\hp{=}+\dot\F F^{kl}\riema \alpha\beta\gamma\delta\n^\alpha x_k^\beta\n^\gamma x_l^\delta h_i^j-\dot\F F
\riema \alpha\beta\gamma\delta\n^\alpha x_i^\beta\n^\gamma x_m^\delta g^{mj}\\
&\hp{=}+ (\F-\tilde f)\riema \alpha\beta\gamma\delta\n^\alpha x_i^\beta\n^\gamma x_m^\delta g^{mj}\\
&\hp{=}+\dot\F F^{kl}\bar R_{\alpha\beta\gamma\delta;\e}\{\n^\alpha x_k^\beta x_l^\gamma x_i^\delta
x_m^\e g^{mj}+\n^\alpha x_i^\beta x_k^\gamma x_m^\delta x_l^\e g^{mj}\}.
\end{aligned}
\end{equation}
\el

Let $\h$ be the vector field \re{h1.18}, or any vector field satisfying \re{h1.17}, and set
\begin{equation}
\tilde v=\spd\h\nu,
\end{equation}
then we have:
\bl[Evolution of $\tilde v$]\lal{h1.4.4}
The quantity $\tilde v$ satisfies the evolution equation
\begin{equation}\lae{h1.4.23}
\begin{aligned}
\dot{\tilde v}-\dot\F F^{ij}\tilde v_{ij}=&\dot\F F^{ij}h_{ik}h_j^k\tilde v
-[(\F-\tilde f)-\dot\F F]\h_{\alpha\beta}\n^\alpha\n^\beta\\
&-2\dot\F F^{ij}h_j^k x_i^\alpha x_k^\beta \h_{\alpha\beta}-\dot\F F^{ij}\h_{\alpha\beta\gamma}x_i^\beta
x_j^\gamma\n^\alpha\\
&-\dot\F F^{ij}\riema \alpha\beta\gamma\delta\n^\alpha x_i^\beta x_k^\gamma x_j^\delta\h_\e x_l^\e g^{kl}\\
&-\tilde f_\beta x_i^\beta x_k^\alpha \h_\alpha g^{ik}.
\end{aligned}
\end{equation}
\el

The derivation is elementary, see the proof of the corresponding lemma in the Lorentzian case, \frl{1.4.4}.

Notice that $\tilde v$ is supposed to satisfy \re{h1.17}, hence
\begin{equation}
\f=-\log(\tilde v-\theta)
\end{equation}
is well defined and there holds
\begin{equation}\lae{h2.5}
\begin{aligned}
\dot \f-\dot\F F^{ij}\f_{ij}=-\{\dot {\tilde v}-\dot\F F^{ij}\tilde v_{ij}\}\frac1{\tilde v-\theta}-\dot\F F^{ij}\f_i\f_j.
\end{aligned}
\end{equation}

Finally, let $\chi$ be the strictly convex function. Its evolution equation is
\begin{equation}\lae{h2.6}
\begin{aligned}
\dot\chi-\dot\F F^{ij}\chi_{ij}&=-[(\F-\tilde f)-\dot \F F]\chi_\al\nu^\al -\dot\F F^{ij}\chi_{\al\bet}x^\al_ix^\bet_j\\
&\le -[(\F-\tilde f)-\dot \F F]\chi_\al\nu^\al -c_0 \dot\F F^{ij}g_{ij}
\end{aligned}
\end{equation}
where $c_0>0$ is independent of $t$.

We can now prove \rt{h1.4}:
\bp[Proof of \rt{h1.4}]
Let $\zeta$ and $w$ be respectively defined by
\begin{align}
\zeta&=\sup\set{{h_{ij}\h^i\h^j}}{{\norm\h=1}},\\
w&=\log\zeta+ \f+\lam \chi,\lae{h2.3.14}
\end{align}
where $\lam>0$ is supposed to be large.  We claim that
$w$ is bounded, if $\lam$ is chosen sufficiently large.

Let $0<T<T^*$, and $x_0=x_0(t_0)$, with $ 0<t_0\le T$, be a point in $M(t_0)$ such
that
\begin{equation}
\sup_{M_0}w<\sup\set {\sup_{M(t)} w}{0<t\le T}=w(x_0).
\end{equation}

We then introduce a Riemannian normal coordinate system $(\x^i)$ at $x_0\in
M(t_0)$ such that at $x_0=x(t_0,\x_0)$ we have
\begin{equation}
g_{ij}=\delta_{ij}\q \tup{and}\q \zeta=h_n^n.
\end{equation}

Let $\tilde \h=(\tilde \h^i)$ be the contravariant vector field defined by
\begin{equation}
\tilde \h=(0,\dotsc,0,1),
\end{equation}
and set
\begin{equation}
\tilde \zeta=\frac{h_{ij}\tilde \h^i\tilde \h^j}{g_{ij}\tilde \h^i\tilde \h^j}\raise 2pt
\hbox{.}
\end{equation}

$\tilde\zeta$ is well defined in neighbourhood of $(t_0,\x_0)$.

Now, define $\tilde w$ by replacing $\zeta$ by $\tilde \zeta$ in \re{h2.3.14}; then, $\tilde w$
assumes its maximum at $(t_0,\x_0)$. Moreover, at $(t_0,\x_0)$ we have 
\begin{equation}
\dot{\tilde \zeta}=\dot h_n^n,
\end{equation}
and the spatial derivatives do also coincide; in short, at $(t_0,\x_0)$ $\tilde \zeta$
satisfies the same differential equation \re{h1.4.1} as $h_n^n$. For the sake of
greater clarity, let us therefore treat $h_n^n$ like a scalar and pretend that $w$
is defined by 
\begin{equation}
w=\log h_n^n+ \f+\lam \chi.
\end{equation} 

From the equations \re{h1.4.1}, \re{h2.5},  \re{h2.6} and \re{h1.8}, we infer, by observing the special form of $\F$, i.e., $\F(F)=F$, $\dot\F=1$, $\tilde f=f$ and using the monotonicity and homogeneity of $F$
\begin{equation}
F=F(\ka_i)=F(\tfrac{\ka_1}{\ka_n},\ldots,1)\ka_n\le F(1,\dots,1)\ka_n
\end{equation}
that in $(t_0,\xi_0)$
\begin{equation}\lae{h2.15}
\begin{aligned}
0&\le-\tfrac12\dot\F F^{ij}h_{ki}h^k_k\frac\theta{\tilde v-\theta}-f h^n_n+c(\theta)\dot\F F^{ij}g_{ij}+\lam c\\
&\hp{\le}\;-\lam c_0\dot\F F^{ij}g_{ij}-\dot \F F^{ij}\f_i\f_j+\dot\F F^{ij}(\log h^n_n)_i(\log h^n_n)_j\\
&\hp{\le}\;+\frac2{\ka_n-\ka_1}\dot\F \sum_{i=1}^n(F_n-F_i)(h_{ni;}^{\hp{ni;}n})^2 (h^n_n)^{-1}.
\end{aligned}
\end{equation}
Similarly as in \cite[p. 197]{cg:scalar}, we distinguish two cases

\cvm
\tit{Case} $1$.\q Suppose that
\begin{equation}
\abs{\ka_1}\ge \e_1 \ka_n,
\end{equation}
where $\e_1>0$ is small, notice that the principal curvatures are labelled according to \re{h1.5}. Then, we infer from \cite[Lemma 8.3]{cg:scalar}
\begin{equation}
F^{ij}h_{ki}h^k_j\ge \tfrac1n F^{ij}g_{ij}\e_1^2\ka_n^2,
\end{equation}
and 
\begin{equation}
F^{ij}g_{ij}\ge F(1,\ldots,1),
\end{equation}
for a proof see e.g., \cite[Lemma 2.2.19]{cg:cp}.

Since $Dw=0$,
\begin{equation}
D\log h^n_n=-D\f-\lam D\chi,
\end{equation}
we obtain
\begin{equation}
\dot\F F^{ij}(\log h^n_n)_i(\log h^n_n)_j=\dot \F F^{ij}\f_i\f_j+2\lam \dot\F F^{ij}\f_i\chi_j+\lam^2\dot\F F^{ij}\chi_i\chi_j,
\end{equation}
where
\begin{equation}
\abs{\f_i}\le c\abs{\ka_i}+c,
\end{equation}
as one easily checks. 

Hence, we conclude that $\ka_n$ is a priori bounded in this case.

\cvm
\tit{Case} $2$.\q Suppose that
\begin{equation}\lae{h2.22}
\ka_1\ge -\e_1\ka_n,
\end{equation}
then, the last term in inequality \re{h2.15} is estimated from above by
\begin{equation}
\begin{aligned}
&\frac2{1+\e_1}\dot\F \sum_{i=1}^n(F_n-F_i)(h_{ni;}^{\hp{ni;}n})^2 (h^n_n)^{-2}&\le  \\
&\frac2{1+2\e_1}\dot\F \sum_{i=1}^n(F_n-F_i)(h_{nn;}^{\hp{nn;}i})^2 (h^n_n)^{-2}\\
&\qq +c(\e_1)\dot\F \sum_{i=1}^{n-1}(F_i-F_n)\ka_n^{-2}
\end{aligned}
\end{equation} 
where we used the Codazzi equation. The last sum can be easily balanced.

The terms in \re{h2.15} containing the derivative of $h^n_n$ can therefore be estimated from above by
\begin{equation}
\begin{aligned}
&-\frac{1-2\e_1}{1+2\e_1}\dot\F \sum_{i=1}^nF_i(h_{nn;}^{\hp{nn;}i})^2 (h^n_n)^{-2}\\
&+\frac2{1+2\e_1}\dot\F F_n\sum_{i=1}^n(h_{nn;}^{\hp{nn;}i})^2 (h^n_n)^{-2}\\
&\le\dot\F F_n\sum_{i=1}^n(h_{nn;}^{\hp{nn;}i})^2 (h^n_n)^{-2}\\
&=\dot\F F_n \norm{D\f+\lam D\chi}^2\\
&=\dot\F F_n\{\norm{D\f}^2+\lam^2\norm{D\chi}^2+2\lam \spd{D\f}{D\chi}\}.
\end{aligned}
\end{equation}

Hence we finally deduce
\begin{equation}
\begin{aligned}
0\le -\dot\F \tfrac12 F_n\ka_n^2\frac\theta{\tilde v-\theta}&+c\lam^2\dot\F F_n(1+\ka_n)-f\ka_n+\lam c \\
&+(c(\theta)-\lam c_0)\dot\F F^{ij}g_{ij}
\end{aligned}
\end{equation}
Thus, we obtain an a priori estimate 
\begin{equation}\lae{h2.28}
\ka_n\le \const,
\end{equation}
if $\lam$ is chosen large enough. Notice that $\e_1$ is only subject to the requirement $0<\e_1<\frac12$.
\ep

\br
Since  the initial condition $F\ge f$ is preserved under the flow, a simple application of the maximum principle, \cf \cite[Lemma 5.2]{cg96}, we conclude that the principal curvatures of the flow hypersurfaces stay in a compact subset of $\C$.
\er

\br\lar{h2.3}
These a priori estimates are of course also valid, if $M$ is a stationary solution.
\er

\subsection{Proof of \rt{h1.5}}

We  consider the curvature flow \re{h1.13} with initial hypersurface $M_0=M_2$. The flow will exist in a maximal time interval $[0,T^*)$ and will stay in $\bar\Om$. We shall also assume that $M_2$ is not already a solution of the problem for otherwise the flow will be stationary from the beginning.

Furthermore, the flow hypersurfaces can be written as graphs
\begin{equation}
M(t)=\graph u(t,\cdot)
\end{equation}
over $\so$, since the initial hypersurface has this property and all flow hypersurfaces are supposed to be convex, i.e., uniform $C^1$-estimates are guaranteed, \cf \cite{cg96}.

The curvature estimates from \rt{h1.4} ensure that the curvature operator is uniformly elliptic, and in view of  well-known regularity results we then conclude that the flow exists for all time and converges in $C^{m+2,\bet}(\so)$ for some $0<\bet\le\al$ to a limit hypersurface $M$, that will be a stationary solution, \cf \cite[Section 6]{cg:survey} and also \frt{6.5.5}.

\section{Foliation of a spacetime by CMC hypersurfaces}

Hypersurfaces of prescribed mean curvature especially those with constant
mean curvature play an important role in general relativity. In \ci{cg1} the
existence of closed hypersurfaces of prescribed mean curvature in a globally
hyperbolic Lorentz manifold with a compact Cauchy hypersurface was proved
provided there were barriers. The proof consisted of two parts, the a priori
estimates for the gradient and the application of a fixed point theorem. That
latter part of the proof was rather complicated.

Ecker and Huisken, therefore, gave another existence proof using an
evolutionary approach, but they had to assume that the timelike convergence
condition is satisfied, and, even more important, that the prescribed mean
curvature satisfies a structural monotonicity condition, cf. \ci{eh1}. These are
serious restrictions which had to be assumed because the authors relied on the
gradient estimate of Bartnik \ci{br:mean}, who had proved another a priori estimate
in the elliptic case.

We later gave an existence proof, using a curvature flow method, that works in an arbitrary globally hyperbolic spacetime without any assumptions on the ambient curvature as long as there are barriers, \cf \cite{cg:mz}.

Let $N$ be a globally hyperbolic Lorentzian manifold with a compact Cauchy
hypersurface $\mc S_0$ and a sufficiently smooth proper time function $x^0$. Consider the problem of finding a closed hypersurface
of prescribed mean curvature $H$ in $N$, or more precisely, let $\Om$ be a
connected open subset of $N$, $f\in C^{0,\al}(\bar \Om)$, then we look for a
hypersurface $M\su \Om$ such that
\begin{equation}
\fv HM=f(x)\qq \A \,x\in M,
\end{equation}
where $\fv HM$ means that  $H$ is evaluated at the vector $(\ka_i(x))$ the
components of which are the principal curvatures of $M$.

We assume that $\pa \Om$ consists of two compact, connected,
spacelike hypersurfaces $M_1$ and $M_2$, where $M_1$ is supposed to lie in the
\tit{past} of $M_2$. The $M_i$ should act as barriers for $(H,f)$, where $M_2$ is an upper and $M_1$ a lower barrier.

Notice that we do not assume $f$ to be positive, hence the mean curvature function is supposed to be defined in $\R[n]$ and not in the usual cone $\C_1$, see  \cite[Definition 1.2.10]{cg:cp}.

In \ci[Section 6]{cg1} we proved the following theorem:

\bt\lat{3.2.2.2}
Let $M_1$ be a lower and $M_2$ be an upper barrier for $(H,f)$, $f\in
C^{0,\al}(\bar \Om)$. Then, the problem 
\begin{equation}\lae{3.2.2.4}
\fv HM=f
\end{equation}
has a solution $M\su \bar \Om$ of class $C^{2,\al}$ that can be written as a graph
over the Cauchy hypersurface $\mc S_0$.
\et

The crucial point in the proof is an a priori estimate in the $C^1$-norm and for
this estimate only the boundedness of $f$ is needed, i.e., even for merely
bounded $f$ $H^{2,p}$-solutions exist.

We want to give a proof of \rt{3.2.2.2} that is based on the curvature flow method, and
to make this method work, we have to assume temporarily slightly higher
degrees of regularity for the barriers and right-hand side, i.e., we assume the
barriers to be of class $C^{6,\al}$ and $f$ to be of class $C^{4,\al}$. We can achieve
these assumptions by approximation without sacrificing the barrier conditions,
cf. \cite[Remark {3.5.2}]{cg:cp}.

To solve \re{3.2.2.4} we look at the evolution problem 
\begin{equation}\lae{3.2.2.5}
\begin{aligned}
\dot x&=(H-f)\n,\\
x(0)&=x_0,
\end{aligned}
\end{equation}
where $x_0$ is an embedding of an initial hypersurface $M_0$, for which we
choose $M_0=M_2$,  $H$ is the mean curvature of the flow hypersurfaces
$M(t)$ with respect to the past directed normal $\n$, and $x(t)$ is an embedding
of $M(t)$, \cf \frt{5.7.0.0}.

The existence result in \rt{3.2.2.2} can be used to prove that a spacetime $N$, satisfying the assumptions of the previous sections, can be foliated by \tit{constant mean curvature}\index{constant mean curvature (CMC)} hypersurfaces, abbreviated (CMC) hypersurfaces, or that at least important parts of $N$, like a future or past end, can be foliated by \ind{CMC hypersurface}s, and that in those parts, the mean curvature of the leaves of the foliation can be used as new smooth time function. 

Of course $N$ has to satisfy some additional conditions in order that the existence of such a foliation can be proved.

If the \tit{\ind{timelike convergence condition}} holds in in $N$, i.e., if
\begin{equation}\lae{3.6.1}
\bar R_{\al\bet}\nu^\al\nu^\bet\ge 0\qq\A\,\spd\nu\nu=-1,
\end{equation}
and if $N$ has future and past mean curvature barriers, see \frd{3.6.1.2} for details, then we proved in \cite{cg1} that $N$ can be foliated by CMC hypersurfaces. The mean curvature of the leaves can then be used as a smooth time function at least in those parts, where the mean curvature of the slices does not vanish, \cf \cite{cg:foliation}.

We later generalized this result by replacing the condition \re{3.6.1} by the weaker assumptions
\begin{equation}\lae{3.6.2}
\bar R_{\al\bet}\nu^\al\nu^\bet\ge -\Lam\qq\A\,\spd\nu\nu=-1,
\end{equation}
where $\Lam\ge 0$ is a constant, and showed that the former results were still valid in future and past ends of $N$, \cf \cite{cg:foliation2}.

We shall first present the foliation results for a spacetime satisfying the preceding weak condition on the Ricci tensor. Setting $\Lam=0$, we then immediately obtain the corresponding results for spacetimes satisfying the timelike convergence condition in those parts of $N$ that are foliated by slices with non-zero CMC hypersurfaces. Only the possible presence of \tit{\ind{maximal hypersurface}s} will require some additional arguments.

Thus let $N$ be a $(n+1)$-dimensional spacetime with a compact Cauchy hypersurface,
so that $N$ is topologically a product, $N=I\times \so$, where $\so$ is a compact
Riemannian manifold and $I=(a,b)$ an interval. 

\bd
A \tit{\ind{future end}} of $N$, in symbols \inds{$N_+$}, is defined by
\begin{equation}
N_+=(x^0)^{-1}[a_0,b)
\end{equation}
and similarly a \tit{\ind{past end}} by
\begin{equation}
N_-=(x^0)^{-1}(a,b_0],
\end{equation}
where $a_0$ and $b_0$ belong to $I$.
\ed

To apply the existence result in \frt{3.2.2.2}, we need barriers, or more precisely, a future (past)  mean curvature barrier, \cf \frd{3.6.1.2}.

Our first results are described in the following two theorems.

\bt\lat{3.6.0.1}
Suppose that in a future end $N_+$ of $N$ the Ricci tensor satisfies the estimate
\re{3.6.2}, and suppose that a future mean curvature barrier exists,  then a
slightly smaller future end $\tilde N_+$ can be foliated by CMC spacelike
hypersurfaces, and there exists a smooth time function $x^0$ such that the slices
\begin{equation}
M_\tau=\{x^0=\tau\},\qq \tau_0<\tau<\un,
\end{equation}
have mean curvature $\tau$ for some $\tau_0>\sqrt{n\Lam}$. The precise value of
$\tau_0$ depends on the mean curvature of a lower barrier.
\et

\bt\lat{3.6.0.2}
Suppose that a future end $N_+=(x^0)^{-1}[a_0,b)$ of $N$ can be covered by a time function
$x^0$ such that the mean curvature of the slices $M_t=\{x^0=t\}$ is
non-negative and the volume of $M_t$ decays to zero
\begin{equation}
\lim_{t\ra b}\abs{M_t}=0,
\end{equation}
then the volume $\abs{M_k}$ of any sequence of spacelike hypersurfaces $M_k$
that approach $b$, i.e.,
\begin{equation}
\lim_k\inf_{M_k}x^0=b,
\end{equation}
decays to zero. Thus, in case the additional conditions of \rt{3.6.0.1} are also satisfied,
the volume of the CMC hypersurfaces $M_\tau$ converges to zero
\begin{equation}
\lim_{\tau\ra\un}\abs{M_\tau}=0.
\end{equation}
\et

$N$ is also future timelike incomplete, if there is a compact spacelike hypersurface
$M$ with mean curvature $H$ satisfying
\begin{equation}
H\ge H_0>\sqrt{n\Lam},
\end{equation}
due to a result in \cite{galloway:cft}.

\subsection{Foliation of future ends}

Let us recall the results in \rb{1.3.6} and \frn{1.3.7}, which, in the present situation, can be phrased  like this:  In a given Gaussian coordinate
system $(x^\al)$ the coordinate slices
$M(t)=\{x^0=t\}$ can be looked at as a solution of the evolution problem
\begin{equation}\lae{3.6.2.1}
\dot x=-e^\psi \nu,
\end{equation}
where $\nu=(\nu^\al)$ is the past directed normal vector. The embedding
$x=x(t,\xi)$ is then given by $x=(t,x^i)$, where $(x^i)$ are local coordinates for
$\so$.

Let $\bar g_{ij},\bar h_{ij}$ and $\bar H$ be the induced metric,
second fundamental and mean curvature of  the coordinate slices,  then
the evolution equations
\begin{equation}\lae{3.6.2.2}
\dot{\bar g}_{ij}=-2 e^\psi \bar h_{ij}
\end{equation}
and
\begin{equation}\lae{3.6.2.3}
\dot{\bar H}=-\D e^\psi +(\abs{\bar A}^2+\bar R_{\al\bet}\nu^\al\nu^\bet) e^\psi
\end{equation}
are valid.

Now, let $M_0$ be a smooth connected spacelike hypersurface and consider in a tubular
neighbourhood $\mc U$ of $M_0$ hypersurfaces $M$ that can be written as graphs
over $M_0$, $M=\graph u$, in the corresponding normal Gaussian coordinate system.
Then the mean curvature of $M$ can be expressed as
\begin{equation}\lae{3.6.2.4}
H=-\D u+\bar H+v^{-2}u^iu^j\bar h_{ij},
\end{equation}
\cf equation \fre{2.16}, and hence, choosing $u=\e\f$, $\f\in C^2(M_0)$, we deduce
\begin{equation}\lae{3.6.2.5}
\begin{aligned}
\frac d{d\e}\fv H {\e=0}&=-\D \f +\dot{\bar H}\f\\[\cma]
&= -\D\f +(\abs{\bar A}^2+\bar R_{\al\bet}\nu^\al\nu^\bet)\f.
\end{aligned}
\end{equation}

Next we shall prove that CMC hypersurfaces are monotonically ordered, if the mean
curvatures are sufficiently large.

\bl\lal{3.6.2.1}
Let $M_i=\graph u_i$, $i=1,2$, be two spacelike hypersurfaces such that the \resp
mean curvatures $H_i$ satisfy
\begin{equation}\lae{3.6.2.6b}
H_1<H_2
\end{equation}
where $H_2$ is constant,\footnote{It would suffice to require $H_1<\inf_{M_2}H_2$.} $H_2=\tau_2$, and
\begin{equation}\lae{3.6.2.7b}
\sqrt{n\Lam}<\tau_2,
\end{equation}
 then there holds
\begin{equation}\lae{3.6.2.7}
u_1<u_2.
\end{equation}
\el

\bp
We first observe that the weaker conclusion
\begin{equation}\lae{3.6.2.10}
u_1\le u_2
\end{equation}
is as good as the strict inequality in \re{3.6.2.7}, in view of the maximum principle.

Hence, suppose that \re{3.6.2.10} is not valid, so that
\begin{equation}\lae{3.6.2.11}
E(u_1)=\set{x\in\so}{u_2(x)<u_1(x)}\ne\eS.
\end{equation}

Then there exist points $p_i\in M_i$ such that
\begin{equation}
0<d_0=d(M_2,M_1)= d(p_2,p_1) =\sup \set{d(p,q)}{(p,q)\in M_2\times M_1},
\end{equation}
where $d$ is the Lorentzian distance function. Let $\f$ be a maximal geodesic from
$M_2$ to $M_1$  realizing this distance with endpoints $p_2$ and $p_1$, and
parametrized by arc length.

Denote by $\bar d$ the Lorentzian distance function to $M_2$, i.e., for $p\in
I^+(M_2)$
\begin{equation}
\bar d(p)=\sup_{q\in M_2}d(q,p).
\end{equation}

Since $\f$ is maximal, $\C=\set{\f(t)}{0\le t<d_0}$ contains no focal points of
$M_2$,
\cf \cite[Theorem 34, p. 285]{bn}, hence there exists an open neighbourhood $\mc
V=\mc V(\C)$ such that $\bar d$ is smooth in $\mc V$, \cf \cite[Theorem {1.9.15}]{cg:cp}. $\mc V$ is part of the largest tubular neighbourhood of $M_2$, and hence covered by an associated normal Gaussian coordinate system $(x^\al)$ satisfying $x^0=\bar d$ in $\{x^0>0\}$, see \cite[Theorem {1.9.22}]{cg:cp}.

Now, $M_2$ is the level set $\{\bar d=0\}$, and the level sets 
\begin{equation}
M(t)=\set{p\in \mc V}{\bar d(p)=t}
\end{equation}
are  smooth hypersurfaces.

Thus, the mean curvature $\bar H(t)$ of $M(t)$ satisfies the equation
\begin{equation}
\dot {\bar H}=\abs{\bar A}^2+\bar R_{\al\bet}\nu^\al\nu^\bet,
\end{equation}
\cf \re{3.6.2.3}, and therefore we have
\begin{equation}\lae{3.6.2.15b}
\dot{\bar H}\ge \tfrac1n\abs{\bar H}^2-\Lam >0,
\end{equation}
in view of \re{3.6.2.7b}.

Next, consider a tubular neighbourhood $\mc U$ of $M_1$ with corresponding
normal Gaussian coordinates $(x^\al)$. The level sets
\begin{equation}
\tilde M(s)=\{x^0=s\},\qq-\e<s<0,
\end{equation}
lie in the past of $M_1=\tilde M(0)$ and are smooth for small $\e$.

Since the geodesic $\f$ is normal to $M_1$, it is also normal to $\tilde M(s)$ and
the length of the geodesic segment of $\f$ from $\tilde M(s)$ to $M_1$ is exactly
$-s$, i.e., equal to the distance from $\tilde M(s)$ to $M_1$, hence we deduce
\begin{equation}
d(M_2,\tilde M(s))=d_0+s,
\end{equation}
i.e., $\set{\f(t)}{0\le t\le d_0+s}$ is also a maximal geodesic from $M_2$ to $\tilde
M(s)$, and we conclude further that, for fixed $s$, the hypersurface $\tilde
M(s)\ii\mc V$ is contained in the past of $M(d_0+s)$ and touches $M(d_0+s)$ in
$p_s=\f(d_0+s)$. The maximum principle then implies
\begin{equation}
\fv H{\tilde M(s)}(p_s)\ge \fv H {M(d_0+s)}(p_s)>\tau_2,
\end{equation}
in view of \re{3.6.2.15b}.

On the other hand, the mean curvature of $\tilde M(s)$ converges to the mean
curvature of $M_1$, if $s$ tends to zero, hence we conclude
\begin{equation}
H_1(\f(d_0))\ge \tau_2,
\end{equation}
contradicting \re{3.6.2.6b}.
\ep

\bc\lac{3.6.2.2}
The CMC hypersurfaces with mean curvature
\begin{equation}
\tau>\sqrt{n\Lam}
\end{equation}
are uniquely determined.
\ec

\bp
Let $M_i=\graph u_i$, $i=1,2$, be two hypersurfaces with mean curvature $\tau$
and suppose, e.g., that
\begin{equation}
\set{x\in\so}{u_1(x)<u_2(x)}\ne\eS.
\end{equation}
Consider a tubular neighbourhood of $M_1$ with a corresponding future oriented
normal Gaussian coordinate system $(x^\al)$. Then the evolution of the mean
curvature of the coordinate slices satisfies
\begin{equation}
\dot{\bar H}=\abs{\bar A}^2+\bar R_{\al\bet}\nu^\al\nu^\bet \ge \frac
1n\abs{\bar H}^2-\Lam >0
\end{equation}
in a neighbourhood of $M_1$, i.e., the coordinate slices $M(t)=\{x^0=t\}$, with
$t>0$, have all mean curvature $\bar H(t)>\tau$. Using now $M_1$ and $M(t)$,
$t>0$, as barriers, we infer from \frt{3.2.2.2} that for any $\tau'\in\R[]$, $\tau<\tau' <\bar H(t)$,
there exists a spacelike hypersurface $M_{\tau'}$ with mean curvature $\tau'$, such
that
$M_{\tau'}$ can be expressed as a graph over $M_1$, $M_{\tau'}=\graph u$, where
\begin{equation}
0<u<t.
\end{equation}

Writing $M_{\tau'}$ as graph over $\so$ in the original coordinate system without
changing the notation for $u$, we obtain
\begin{equation}\lae{3.6.2.14}
u_1<u,
\end{equation}
and, by choosing $t$ small enough, we may also conclude that
\begin{equation}\lae{3.6.2.15}
E(u)=\set{x\in\so}{u(x)<u_2(x)}\ne \eS,
\end{equation}
which is impossible, in view of the preceding result.
\ep

\bl\lal{3.7.3}
Under the assumptions of \rt{3.6.0.1}, let $M_{\tau_0}=\graph u_{\tau_0}$ be a CMC
hypersurface with mean curvature $\tau_0>\sqrt{n\Lam}$, then the future of
$M_{\tau_0}$ can be foliated by CMC hypersurfaces
\begin{equation}\lae{3.6.2.28}
I^+(M_{\tau_0})=\uuu_{\tau_0<\tau<\un}M_\tau.
\end{equation}
The $M_\tau$ can be written as graphs over $\so$
\begin{equation}
M_\tau=\graph u(\tau,\cdot),
\end{equation}
such that $u$ is strictly monotone increasing with respect to $\tau$, and continuous
in $[\tau_0,\un)\times \so$.
\el

\bp
The monotonicity and continuity of $u$ follows from \rl{3.6.2.1} and \rc{3.6.2.2}, in view of
the a priori estimates.

Thus, it remains to verify the relation \re{3.6.2.28}. Let $p=(t,y^i)\in I^+(M_{\tau_0})$,
then we have to show $p\in M_\tau$ for some $\tau>\tau_0$.

From the existence result in  \rt{3.2.2.2}  we deduce  that there exists a family of CMC
hypersurfaces $M_\tau$
\begin{equation}
\set{M_\tau}{\tau_0\le \tau<\un},
\end{equation}
since there is a future mean curvature barrier.

Define $u(\tau,\cdot)$ by
\begin{equation}
M_\tau=\graph u(\tau,\cdot),
\end{equation}
then we have
\begin{equation}
u(\tau_0,y)<t<u(\tau^*,y)
\end{equation}
for some large $\tau^*$, because of the mean curvature barrier condition, which,
together with \rl{3.6.2.1}, implies that the CMC hypersurfaces run into the future
singularity, if $\tau$ goes to infinity.

In view of the continuity of $u(\cdot,y)$ we conclude that there exists $\tau_1$ such that $\tau_0<\tau_1<\tau^*$ and 
\begin{equation}
u(\tau_1,y)=t,
\end{equation}
 hence
$p\in M_{\tau_1}$.
\ep

\br
The continuity and monotonicity of $u$ holds in any coordinate system $(x^\al)$,
even in those that do not cover the future completely like the normal Gaussian
coordinates associated with a spacelike hypersurface, which are defined in a tubular
neighbourhood.
\er

The proof of \frt{3.6.0.1} is now almost finished. The remaining arguments are given in several steps.

 We have to show that the mean curvature parameter $\tau$ can be
used as a time function in $\{\tau_0<\tau<\un\}$, i.e., $\tau$ should be smooth with
a non-vanishing gradient. Both properties are local properties.

\hinweis{First step} \lan{3.7.5}
 Fix an arbitrary $\tau'\in (\tau_0,\un)$, and consider a tubular
neighbourhood $\mc U$ of $M'=M_{\tau'}$. The $M_\tau\su \mc U$ can then be
written as graphs over $M'$, $M_\tau=\graph u(\tau,\cdot)$. For small $\e>0$ we
have
\begin{equation}
M_\tau\su \mc U\qq\A\,\tau\in (\tau'-\e,\tau'+\e)
\end{equation}
and with the help of the implicit function theorem we shall show that $u$ is smooth.
Indeed, define the operator $G$
\begin{equation}
G(\tau,\f)=H(\f)-\tau,
\end{equation}
where $H(\f)$ is an abbreviation for the mean curvature of $\graph \fv \f{M'}$. Then
$G$ is smooth and from \re{3.6.2.5} we deduce that $D_2G(\tau',0)\f$ equals
\begin{equation}\lae{3.7.37}
-\D\f+(\norm{ A}^2+\bar R_{\al\bet}\nu^\al\nu^\bet)\f,
\end{equation}
where the Laplacian, the second fundamental form and the normal correspond to $M'$.
Hence $D_2G(\tau',0)$ is an isomorphism and the implicit function theorem implies
that $u$ is smooth.

\hinweis{Second step}\lan{3.7.6}
 Still in the tubular neighbourhood of $M'$, define the coordinate
transformation
\begin{equation}
\F(\tau,x^i)=(u(\tau,x^i),x^i);
\end{equation}
note that $x^0=u(\tau,x^i)$. Then we have
\begin{equation}
\det D\F=\frac{\pa u}{\pa \tau}=\dot u.
\end{equation}

$\dot u$ is non-negative; if it were strictly positive, then $\F$ would be a
diffeomorphism, and hence $\tau$ would be smooth with non-vanishing gradient. To
prove $\dot u>0$, observe that the CMC hypersurfaces in $\mc U$ satisfy an equation
\begin{equation}
H(u)=\tau,
\end{equation}
where the left hand-side can be expressed as in \re{3.6.2.4}. Differentiating both sides
with respect to $\tau$ and evaluating for $\tau=\tau'$, i.e., on $M'$, where
$u(\tau',\cdot)=0$, we get
\begin{equation}\lae{3.7.1}
-\D\dot u+(\abs{A}^2+\bar R_{\al\bet}\nu^\al\nu^\bet)\dot u=1.
\end{equation}

In a point, where $\dot u$ attains its minimum, the maximum principle implies
\begin{equation}
(\abs{A}^2+\bar R_{\al\bet}\nu^\al\nu^\bet)\dot u\ge1,
\end{equation}
hence $\dot u\ne0$ and $\dot u$ is therefore strictly positive.

\br
The results in \frt{3.6.0.1} are also valid in a past end, if $N$ has a past mean curvature barrier. Moreover, the assumption in the future (past) mean barrier condition that the mean curvature of the barriers converge to $\un$ \resp $-\un$ can be easily replaced by the assumption that the limits are finite numbers as long as the absolute values of these numbers are strictly larger than $\sqrt{n\Lam}$.

If $\Lam=0$, the mean curvature of future \resp past barriers are also allowed to converge to $0$.
\er

\subsection{Proof of \rt{3.6.0.2}}

Let $x^0$ be time function satisfying the assumptions of \frt{3.6.0.2}, i.e.,
$N_+=\{a_0<x^0<b\}$, the mean curvature of the slices $M(t)=\{x^0=t\}$ is
non-negative, and
\begin{equation}
\lim_{t\ra b}\abs{M(t)}=0,
\end{equation}
and let $M_k$ be a sequence of connected, spacelike hypersurfaces such that
\begin{equation}
\lim\inf_{M_k}x^0=b.
\end{equation}

Let us write $M_k=\graph u_k$ as graphs over $\so$. Then
\begin{equation}
g_{ij}=e^{2\psi}(u_iu_j+\s_{ij}(u,x))
\end{equation}
is the induced metric, where we dropped the index $k$ for better readability, and the
volume element of $M_k$ has the form
\begin{equation}
d\m=v\sqrt{\det(\bar g_{ij}(u,x))}\,dx,
\end{equation}
where
\begin{equation}\lae{3.6.3.5}
v^2=1-\s^{ij}u_iu_j<1,
\end{equation}
and $(\bar g_{ij}(t,\cdot))$ is the metric of the slices $M(t)$.

From \re{3.6.2.2} we deduce
\begin{equation}\lae{3.6.3.6}
\frac d{dt}\sqrt{\det(\bar g_{ij}(t,\cdot))}=-e^\psi\bar H\sqrt{\det(\bar g_{ij})}\le 0.
\end{equation}

Now, let $a_0<t<b$ be fixed, then for \aev $k$ we have
\begin{equation}\lae{3.6.3.7}
t<u_k
\end{equation}
and hence
\begin{equation}
\begin{aligned}
\abs{M_k}&=\int_\so v \sqrt{\det(\bar g_{ij}(u_k,x))}\,dx\\[\cma]
&\le \int_\so\sqrt{\det(\bar g_{ij}(t,x)}\,dx =\abs{M(t)},
\end{aligned}
\end{equation}
in view of \re{3.6.3.5}, \re{3.6.3.6} and \re{3.6.3.7}, and we conclude
\begin{equation}
\limsup\abs{M_k}\le \abs{M(t)}\qq\A\,a_0<t<b,
\end{equation}
and thus
\begin{equation}
\lim\abs{M_k}=0.
\end{equation}

\subsection{The case $\Lam=0$}

Suppose now that $N$ satisfies the timelike convergence condition and assume that there exist closed, spacelike hypersurfaces with strictly positive and strictly negative mean curvature. Then there exists a real number $\e_0>0$ and a family of $\mc M_{\e_0}$ of closed spacelike graphs $M_\tau$ of mean curvature $\tau$ for any $\tau\in [-\e_0,\e_0]$, in view of the preceding results.

The hypersurfaces can be written as graphs over $\so$, $M_\tau=\graph u(\tau,\cdot)$, and
\begin{equation}
\tau_1<\tau_2\ne 0 \im u(\tau_1)<u(\tau_2),
\end{equation}
in view of \frl{3.6.2.1}.

In view of the a priori estimates in \cite{cg1} or \cite{cg:mz}, \cf also \frt{3.2.2.2},  the preceding monotonicity relation yields that the limit functions 
\begin{equation}
u_1=\lim_{\tau\ua 0}u(\tau)\q\wed\q u_2=\lim_{\tau\da 0}u(\tau)
\end{equation}
are smooth functions the graphs of which are spacelike maximal hypersurfaces. 

Moreover, any other maximal hypersurface $M=\graph u$ must satisfy
\begin{equation}
u_1\le u\le  u_2.
\end{equation}

The second inequality of this relation follows immediately from \frl{3.6.2.1} applied to $u$ and any $u(\tau)$ with $\tau>0$, which in turn also proves the first inequality by switching the light cones.

\bt\lat{3.8.1}
Assume that $u_1\ne u_2$, then both hypersurfaces are totally geodesic and the metric in the region \inds{$\mc C_0$} of $N$ determined by
\begin{equation}
\mc C_0=\set{(x^0,x)}{u_1\le x^0\le u_2}
\end{equation}
is stationary, i.e., the tubular neighbourhood $\mc U$ of $M_1=\graph u_1$ covers $\mc C_0$ and in the corresponding normal Gaussian coordinate system $(x^\al)$ the metric has the form
\begin{equation}
d\bar s^2=-(dx^0)^2+\s_{ij}(x)dx^idx^j,
\end{equation}
where $\s_{ij}$ is the induced metric of $M_1$ and is hence independent of $x^0$. The hypersurface $M_2$ is a level hypersurfaces in the new coordinate system
\begin{equation}
M_2=\{x^0=t_2\},
\end{equation}
and  the slices 
\begin{equation}
M_t=\{x^0=t\}\qq 0\le t\le t_2,
\end{equation}
which foliate $\mc C_0$, are all totally geodesic.

Thus, a foliation of $N$ is given by
\begin{equation}
\mc C_0\uud (M_\tau)_{\tau\ne0},
\end{equation}
where the family $(M_\tau)_{\tau\ne0}$ is the foliation of $N\sminus \mc C_0$ by CMC hypersurfaces with non-vanishing mean curvature, the existence of which has been proved in \rl{3.7.3}.\footnote{Formally, a foliation has only been proved in the future end $0<\tau_0\le \tau<\un$, but it can obviously be extended to cover $0<\tau<\un$, and similarly for the past end.}
\et

\bp
We first note that, in view of the maximum principle, there holds either $u_1<u_2$ or $u_1=u_2$, hence $u_1<u_2$ and their Lorentzian distance $d_0$ is positive.

Consider now a tubular neighbourhood $U_\e$ of $M_1$ for small $\e$, where $\e$ refers to the upper bound of the signed Lorentzian distance from $M_1$, \cf \cite[Theorem 1.3.13]{cg:cp}. We are actually more interested in the future part of $U_\e$, which is denoted by $U_\e^+$ and consists of those points in $U_\e$ which lie in the future of $M_1$.

Thus, we stipulate that in this proof $U_\e$ should be defined as
\begin{equation}
U_\e=U_{\e_1}^-\uu M_1\uu U_\e^+,
\end{equation}
where $\e_1>0$ is fixed and small, and $\e$ is a variable parameter, satisfying 
\begin{equation}
\e_1\le \e<d_0
\end{equation}
which can be chosen as large as $d_0$, as we shall show. 

Let $(x^\al)$ be the normal Gaussian coordinate system associated with the tubular neighbourhood of $M_1$, i.e., $x^0$ denotes the signed Lorentzian distance from $M_1$ and
\begin{equation}
U_\e^+=\set{p\in U_\e}{0<x^0(p)<\e},
\end{equation}
and the metric in $U_\e$ can be expressed as 
\begin{equation}\lae{3.8.12}
d\bar s^2=-(dx^0)^2+\s_{ij}(x^0,x)dx^idx^j.
\end{equation}

Denote the coordinate slices $\{x^0=t\}$, $0\le t<\e$, by $M(t)$, then these slices can also be written as graphs over the Cauchy hypersurface $\so$ in the original coordinate system
\begin{equation}
M(t)=\fv{\graph u(t)}\so.
\end{equation}

Since $M(0)=M_1$ there holds $u(t)<u_2$,
if $0\le t$ is small, and we shall consider only those $\e$ such that 
\begin{equation}\lae{3.8.14}
u(t)<u_2 \qq\A\,0\le t<\e.
\end{equation}

We claim that all slices $M(t)$ contained in $U_\e$ with $t\ge 0$ are totally geodesic and that the metric $\s_{ij}$ in \re{3.8.12} is independent of $x^0$.

To prove this claim, let $\bar g_{ij}$, $\bar h_{ij}$, $\bar H$ and $\nu$ be the corresponding geometric quantities of $M(t)$. The mean curvature satisfies the evolution equation
\begin{equation}
\dot {\bar H}=\abs{\bar A}^2+\bar R_{\al\bet}\nu^\al\nu^\bet,
\end{equation}
\cf \fre{3.6.2.3} and observe that $\psi=0$.

Hence the mean curvature is non-decreasing, i.e., $\bar H(t)\ge 0$. If one of the $M(t)$, say $M(t_0)$, would be not totally geodesic, then the linearization of the mean curvature operator, evaluate at $M(t_0)$ would be an isomorphism, \cf \fre{3.6.2.5}, and the inverse function theorem would yield the existence of  a hypersurface $M=\fv{\graph u}\so$ in a small neighbourhood of $M(t_0)$ such that
\begin{equation}
\fv H M>\bar H(t_0)\ge 0\q\wed\q u<u_2,
\end{equation}
contradicting the results of  \frl{3.6.2.1}; notice that the mean curvature $H_2$ in that lemma need not be constant, it suffices, if the inequality
\begin{equation}
H_1<\inf_{M_2}H_2
\end{equation}
is valid, since this is all that is needed for  the arguments in the proof.

Thus all hypersurfaces $M(t)$ are totally geodesic and hence the metric $\s_{ij}$ independent of $x^0$, because of the evolution equation \fre{3.6.2.2}. In view of the a priori estimates the slices $M(t)$ are uniformly smooth and the tubular neighbourhood $U_\e$ exists for all $\e$ until the inequality \re{3.8.14} is violated, which will only be the case, if $\e> d_0$, for let $\e\le d_0$ and suppose that $0<t_0<\e$ is the first $t$ such that $M(t_0)$ touches $M_2$. Since both hypersurfaces are maximal, the maximum principle would yield $M(t_0)=M_2$, a contradiction, since $t_0<d_0$ and $t_0$ is also the Lorentzian distance of $M(t_0)$ to $M_1$.
\ep

\br
The mean curvature of the CMC leaves $M_\tau$, $\tau\ne 0$, can be used as smooth time function. If $N$ contains just one maximal hypersurface $M_0$, then $\tau$ is smooth in all of $N$ unless $M_0$ is totally geodesic, as can be easily deduced from the arguments in \frn{3.7.5}, where the differential operator in \re{3.7.37} has to be injective, which will be the case, if $M_0$ is not totally geodesic. 
\er

\section{The inverse mean curvature flow in Lorentzian spaces} 
Let us now consider the inverse mean curvature flow (IMCF)
\begin{equation}\lae{7.1}
\dot x=-H^{-1}\nu
\end{equation}
with initial hypersurface $M_0$  in a globally hyperbolic spacetime $N$ with compact Cauchy hypersurface $\so$.

$N$ is supposed to satisfy the timelike convergence condition
\begin{equation}
\bar R_{\al\bet}\nu^\al\nu^\bet\ge0\qq\A\,\spd\nu\nu=-1.
\end{equation}
Spacetimes with compact Cauchy hypersurface that satisfy the timelike convergence condition are also called \tit{cosmological spacetimes}, a terminology due to Bartnik.

In such spacetimes the inverse mean curvature flow will be smooth as long as it stays in a compact set, and, if $\fv H{M_0}>0$ and if the flow exists for all time, it will necessarily run into the future singularity, since the mean curvature of the flow hypersurfaces will become unbounded and the flow will run into the future of $M_0$. Hence the claim follows from \frr{6.4.2}.

However, it might be that the flow will run into the singularity in finite time. To exclude this behaviour we introduced in \cite{cg:imcf} the so-called \tit{strong volume
decay condition}, \cf \rd{5.1.2}. A strong volume decay condition is both
necessary and sufficient in order that the IMCF exists for all time.

\bt\lat{7.1}
Let $N$ be a cosmological spacetime with compact Cauchy hypersurface $\so$ and
with a  future mean curvature barrier. Let $M_0$ be a closed, connected, spacelike hypersurface
with positive mean curvature and assume furthermore that $N$ satisfies a future
volume decay condition. Then the IMCF \re{7.1} with initial hypersurface $M_0$ exists
for all time and provides a foliation of the future \inds{$D^+(M_0)$} of $M_0$.

The evolution parameter $t$ can be chosen as a new time function. The flow
hypersurfaces $M(t)$ are the slices $\{t=\const\}$ and their volume satisfies
\begin{equation}
\abs{M(t)}=\abs{M_0} e^{-t}.
\end{equation}

Defining a new time function $\tau$ by choosing
\begin{equation}
\tau=1-e^{-\frac1n t}
\end{equation}
we obtain  $0\le \tau <1$,
\begin{equation}
\abs{M(\tau)}=\abs{M_0} (1-\tau)^n,
\end{equation}
and  the future singularity corresponds to $\tau=1$.

Moreover, the length $L(\ga)$ of any future directed curve $\ga$ starting from
$M(\tau)$ is bounded from above by
\begin{equation}
L(\ga)\le c (1-\tau),
\end{equation}
where $c=c(n, M_0)$. Thus, the expression $1-\tau$ can be looked at as the radius
of the slices $\{\tau=\const\}$ as well as a measure of the remaining life span of the
spacetime.
\et

Next we shall define the \ind{strong volume decay condition}.

\bd\lad{5.1.2}
Suppose there exists a time function $x^0$ such that the future end of $N$ is
determined by $\{\tau_0\le x^0<b\}$ and the coordinate slices
$M_\tau=\{x^0=\tau\}$ have positive mean curvature with respect to the past
directed normal for $\tau_0\le\tau<b$. In addition the volume $\abs{M_\tau}$
should satisfy
\begin{equation}\lae{5.1.19}
\lim_{\tau\ra b}\abs{M_\tau}=0.
\end{equation}

A decay like that is normally associated with a future singularity and we simply call it
\tit{\ind{volume decay}}. If $(g_{ij})$ is the induced metric of $M_\tau$ and
$g=\det(g_{ij})$, then we have
\begin{equation}\lae{5.1.20}
\log g(\tau_0,x)-\log g(\tau,x)=\int_{\tau_0}^\tau 2 e^\psi \bar H(s,x)\q\A\,x\in \so,
\end{equation}
where $\bar H(\tau,x)$ is the mean curvature of $M_\tau$ in $(\tau,x)$. This relation can be easily derived from the relation \fre{1.3.6} and \frr{1.3.5}.  A detailed proof
is given in \cite{cg:volume}.

In view of \re{5.1.19} the left-hand side of this equation tends to infinity if $\tau$
approaches $b$ for \aev $x\in \so$, i.e.,
\begin{equation}
\lim_{\tau\ra b}\int_{\tau_0}^\tau e^\psi \bar H(s,x)=\un\q \tup{for
\aev}\; x\in\so.
\end{equation}

Assume now, there exists a continuous, positive function $\f=\f(\tau)$ such that
\begin{equation}\lae{5.1.22}
e^\psi \bar H(\tau,x)\ge \f(\tau)\qq\A\, (\tau, x)\in (\tau_0,b)\times \so,
\end{equation}
where
\begin{equation}\lae{5.1.23}
\int_{\tau_0}^b \f(\tau)=\un,
\end{equation}
then we say that the future of $N$ satisfies a \tit{strong volume decay condition}.
\ed

\br
(i) By approximation we may assume that the function $\f$ above is
smooth.

(ii) A similar definition holds for the past of $N$ by simply reversing the time
direction. Notice that in this case the mean curvature of the coordinate slices has to
be negative.
\er

\bl\lal{5.1.4}
Suppose that the future of $N$ satisfies a strong volume decay condition, then there
exist a time function $\tilde x^0=\tilde x^0(x^0)$, where $x^0$ is the time function
in the strong volume decay condition, such that the mean curvature $\bar H$ of the
slices $\tilde x^0=\const$ satisfies the estimate
\begin{equation}\lae{5.1.24}
e^{\tilde\psi}\bar H\ge 1.
\end{equation}
The factor $e^{\tilde\psi}$ is now the conformal factor in the representation
\begin{equation}\lae{5.1.25}
d\bar s^2=e^{2\tilde\psi}(-(d\tilde x^0)^2+\s_{ij}dx^idx^j).
\end{equation}

The range of $\tilde x^0$ is equal to the interval $[0,\un)$, i.e., the singularity
corresponds to $\tilde x^0=\un$.
\el

A proof is given in \cite[Lemma 1.4]{cg:imcf}.

\br
\rt{7.1} can be generalized to spacetimes satisfying
\begin{equation}
\bar R_{\al\bet}\nu^\al\nu^\bet\ge-\Lam\qq\A\,\spd\nu\nu=-1
\end{equation}
with a constant $\Lam\ge 0$, if the mean curvature of the initial hypersurface $M_0$ is sufficiently large
\begin{equation}
\fv H{M_0}>\sqrt{n\Lam},
\end{equation}
\cf \cite{heiko:diplom}. In that thesis it is also shown that the future mean curvature barrier assumption can be dropped, i.e., the strong volume decay condition is sufficient to prove that the IMCF exists for all time and provides a foliation of the future of $M_0$. Hence, the strong volume decay condition already implies the existence of a future mean curvature barrier, since the leaves of the IMCF define  such a barrier.
\er


 \providecommand{\bysame}{\leavevmode\hbox to3em{\hrulefill}\thinspace}
\providecommand{\href}[2]{#2}



\end{document}